\documentclass[12pt,reqno]{amsart}
\usepackage{amscd,amssymb,amsthm,euscript}
\chardef\bslash=`\\ % p. 424, TeXbook
\newtheorem{thm}{Theorem}[section]
\newtheorem{cor}[thm]{Corollary}
\newtheorem{lem}[thm]{Lemma}
\newtheorem{prop}[thm]{Proposition}

\newtheorem*{thm*}{Theorem}
\newtheorem*{lem*}{Lemma}
\newtheorem*{cor*}{Corollary}

\theoremstyle{definition}
\newtheorem{defn}[thm]{Definition}
\newtheorem{ex}[thm]{Example}
\newtheorem{exs}[thm]{Examples}

\theoremstyle{remark}
\newtheorem{rem}[thm]{Remark}

\numberwithin{equation}{section}

\newcommand{\h}{\mathcal{H}}
\newcommand{\ck}{\mathcal{K}}
\newcommand{\ga}{\gamma}
\newcommand{\st}{\sigma}

\def\D{\displaystyle}
\def\ra{\rightarrow}
\def\r{\rangle}
\def\l{\langle}

\font\bb=msbm10 at 11 pt
\font\bbb=msbm10 at 7 pt

\def\croi{\hbox{$\mathop{\mathrel\times\joinrel\mathrel{\vrule height 5pt
 depth 0pt}}$}} 
\def\lcroi{\mbox{$\mathop{\mathrel{\vrule height 5pt
 depth 0pt}\joinrel\mathrel\times} $}}

\def\C{\hbox{\bb C}}
\def\R{\hbox{\bb R}}
\def\N{\hbox{\bb N}}

\def\Z{\hbox{\bb Z}}

\def\n{\hbox{\bbb N}}

\newcommand{\eval}[2][\right]{\relax
  \ifx#1\right\relax \left.\fi#2#1\rvert}

\title{Cohomology of property $T$ groupoids and applications}
\author{Claire Anantharaman-Delaroche}
\address{D\'epartment de Math\'ematiques, Universit\'e d'Orl\'eans,
45067 Orl\'eans, France}
\email{claire@labomath.univ-orleans.fr}
\dedicatory{} 


\begin{document}

\begin{abstract}

We extend the Delorme-Guichardet characterization
of Kazhdan property $T$ groups to $r$-discrete measured groupoids. We
give several applications, in particular to stability
results of Kazhdan property $T$ and to the study of cocycles 
taking  values in a group having the Haagerup property.
                         
\end{abstract} 
\maketitle

\renewcommand{\sectionmark}[1]{}

\section{Introduction}

The notion of Kazhdan property $T$ for discrete group actions and measured
equivalence relations was introduced by Zimmer \cite{zi:cohom}
some twenty years ago, in order to study
the low dimensional cohomology theory of ergodic actions of semi-simple Lie groups
and their lattice subgroups.

This notion extends easily to measured groupoids. However, so
far, the only known examples of measured equivalence relations with property $T$ are
obtained from actions of discrete property $T$ groups. 

Recent striking results of Gaboriau \cite{ga} and Popa \cite{po} bring out
the importance of property $T$ groupoids 
 in ergodic and operator algebras theories. In fact, our paper arises from
a question raised by Damien Gaboriau, namely whether property $T$ of a discrete
measured equivalence relation could be characterized by a fixed point property
relative to each of its actions by affine isometries of affine Hilbert spaces.
We give here a positive answer to this question for $r$-discrete measured groupoids,
more generally.
As we shall see, this  gives a new insight into the subject and, in particular, new
tools to detect property $T$.

The results of this paper are set up for general ergodic $r$-discrete measured
groupoids. However, for the reader's convenience, in this introduction we limit
the description of most of our results to the case of discrete group actions. 
As said before, discrete measured equivalence relations are also 
particular cases of our study. The general situation requires some familiarity with
basic notions on groupoids, that are recalled in the preliminaries.

Throughout this introduction, $(X,\mu)$ will be a standard Borel measure space, and
we consider an ergodic left action of a locally compact group $G$ on $(X,\mu)$,
preserving the class of the measure
$\mu$. In many statements we shall need to assume that $G$ is discrete, to avoid
technical problems. 

A {\it representation} of the associated groupoid
$(X\croi\, G,\mu)$ is a Borel map from $X\croi\, G$ into the unitary group of a
Hilbert space $\ck$, such that for every $g_1,g_2 \in G$,
$$L(x,g_1g_2) = L(x,g_1)L(g_{1}^{-1}x,g_2) \quad \hbox{a.e.}$$
One says that $L$ almost contains unit invariant sections if there exists a net of
Borel maps $\xi_n$ from $X$ into the unit sphere of $\ck$ such that, for every $g\in
G$,
$$\lim_n L(x,g) \xi_n(g^{-1}x) - \xi_n(x) = 0\quad \hbox{a.e.}$$
 Recall (\cite{zi:cohom}, \cite{mo:ergo}) that the $G$-space
$(X,\mu)$ has {\it property} $T$ (or is {\it Kazhdan}) if every representation $L$
almost having invariant unit sections actually has an  invariant unit section, that
is, a Borel map $\xi$ from $X$ into the set of unit vectors in $\ck$ such that for
every
$g\in G$,
$$ L(x,g) \xi(g^{-1}x)= \xi(x)  \quad \hbox{a.e.}$$
(see Section 4.1).

A $L$-{\it cocycle} is a Borel map $b : X\croi\, G \ra \ck$ such that for every
$g_1,g_2 \in G$,
$$b(x,g_1g_2) = b(x,g_1)+ L(x,g_1)b(g_{1}^{-1}x,g_2) \quad \hbox{a.e.},$$
and if there is a Borel map $\xi: X \ra \ck$ such that for $g\in G$,
$$b(x,g) = \xi(x) - L(x,g)\xi(g^{-1}x) \quad \hbox{a.e.},$$
we say that $b$ is a $L$-{\it coboundary}.

If we introduce 
$$\alpha(x,g) : \eta \in \ck \mapsto L(x,g)\eta + b(x,g) \in \ck ,$$
let us note that $\alpha$ is an action of the groupoid $X\croi\, G$ by affine isometries of $\ck$, and that $b$
is the $L$-coboundary defined by $\xi$ if and only if $\xi$ is an invariant section for
$\alpha$, that is, for $g\in G$, one has $\alpha(x,g)\xi(g^{-1}x) = \xi(x)$
a.e. 

The first cohomology group $H^1((X\croi\, G, \mu), L)$ with respect to $L$
is defined in an obvious way. This group gives useful informations on the dynamics.
Our main result is the extension to our setting of the Delorme-Guichardet
characterization of Kazhdan groups, by the
vanishing of the first cohomology group, for every representation
(see \cite[Th\'eor\`eme 7, page 47]{hv}). We follow the approach of Serre, as exposed
in \cite{hv}. An important step is the following characterization of
$L$-coboundaries.

\begin{thm*}[\ref{triv-cocy}] Assume that $G$ is a discrete group. A $L$-cocycle $b$
is a
$L$-coboundary if and only if there exists a subset $E\subset X$ of positive measure
such that for every $x\in E$, we have $\sup\{\|b(x,g)\| : x\in g E\} < +\infty$.
\end{thm*}

From that result, we get:

\begin{thm*}[\ref{cara-suf} and \ref{cara-nec}] Assume that $G$ is a discrete
group. An ergodic $G$-space $(X,\mu)$ has property $T$ if and only if for every
representation $L$ we have $H^1((X\croi\, G, \mu), L)=0$, or equivalently, if and only if
every action of the groupoid $X\croi\, G$ by affine isometries on a Hilbert space
has an invariant section.
\end{thm*}

Section 5 contains several illustrations  of the above characterization. We shall describe
some of them in the rest of this introduction.

Before, observe that, as expected, property $T$ actions are essentially incompatible
with amenability in the sense of Zimmer. Let us spell out this fact. Given an
equivariant map $p: (Z,\tau)\ra (Y,\nu)$ between $G$-spaces we say  that
 $(Z,\tau)$ is a {\it measure preserving extension} of $(Y,\nu)$,
 or that  $\big((Z,\tau), (Y,\nu)\big)$ (often written $(Z,Y)$ for simplicity) is a {\it measure preserving pair},
 if $\tau$ admits the disintegration
$\tau =\int \rho^y d\nu(y)$ where $(\rho^y)$ is a family of probability
measures and $g\rho^y = \rho^{gy}$ almost everywhere.  
   If $\big((Y\times G, \nu\times
\lambda), (Y,\nu)\big)$ is a measure preserving extension for the first projection,
where $\lambda$ is the Haar measure on $G$, we say that 
$(Y,\nu)$ is a {\it proper  $G$-space}. This is the analogue
of proper actions in the topological setting. In fact, every ergodic proper 
  $G$-space is essentially transitive with compact stabilizers, and therefore
it is the left action of $G$ on $G/K$ where $K$ is a compact subgroup of $G$
(see Remark \ref{ess-ergo}). We have:

\begin{thm*}[\ref{proper-am}] The following conditions are equivalent:
\begin{itemize}
\item[(i)] $(X,\mu)$ is a proper  $G$-space.
\item[(ii)] $(X,\mu)$ is both an amenable and a Kazhdan $G$-space.
\end{itemize}
\end{thm*}

Let us assume from now on in this introduction that $G$ is discrete. As a first
consequence of our  Theorems \ref{cara-suf} and
\ref{cara-nec}, we  easily get the invariance of property $T$ by similarity. This
gives another way (Remark \ref{ME}) of
understanding Furman's result stating that, for discrete groups, Kazhdan property is
a measure equivalence invariant \cite[Corollary 1.4]{fu1}. 

A second consequence is the fact that the Mackey range of  an usual cocycle
(called homomorphism in this paper) $\alpha :X\croi\, G \ra H$ is an ergodic
$H$-space with property $T$, when $X$ is a property $T$ ergodic $G$-space and
$H$ is a locally compact group (see Corollary \ref{invMackey}). One has to show that
every affine representation of $X\croi\, G$ by isometries of a Hilbert space has an
invariant section. This new formulation allows to follow the pattern used by Zimmer
in his proof of the important fact that amenability is preserved by Mackey range
(see \cite{zi:amen}). 

In particular, if the $G$-space $X$
has property $T$ and if $H$ is amenable, then the Mackey range of $\alpha$ is a
proper  $H$-space by Theorem \ref{proper-am}, so that it is the action by
translation on $H/K$, where $K$ is a compact subgroup of $H$. Thus we obtain a new
proof of the result, due to Zimmer \cite{zi:cohom} and Schmidt \cite{sch3}, stating
that a cocycle
$\alpha : X\croi\, G \ra H$ is cohomologous to a cocycle into a compact subgroup of
$H$, under the assumptions that the $G$-action has property $T$ and $H$
is amenable. We also improve the result of Adams and Spatzier \cite[Corollary
1.6]{as} showing that when $X$ is a property $T$ ergodic $G$-space, $H$ a locally
compact group and $\alpha$ a cocycle such that the skew-product $X\times_\alpha H$
is an ergodic $G$-space, then $H$ has Serre's property $(FA)$ (see \cite[page 70]{hv}
for the definition). In fact, we show in Corollary \ref{skewergo} that $H$ has
property $T$, which is stronger than property $(FA)$.

Next, we turn to the study of several other stability properties of Kazhdan
groupoids, for the proof of which the following lemma is a key result.

\begin{lem*}[\ref{injpair}] Let $\big((Y,\nu),(X,\mu)\big)$ be  a measure preserving 
pair of ergodic $G$-spaces and denote by $\phi$ the groupoid homomorphism
$(y,g) \mapsto (p(y),g)$ where $p: Y \ra X$ is the underlying $G$-equivariant
map. Then for every representation $L$ of
$(X\croi\, G,\mu)$, the natural map 
$$\phi_* : H^1((X\croi\, G, \mu), L) \ra H^1((Y\croi\, G, \nu), L\circ \phi)$$
is injective.
\end{lem*}

In particular, we have the following consequences.

\begin{thm*}[\ref{quotient}] Let $(Y,X)$ be a measure preserving  pair of
ergodic $G$-spaces. Then $(Y,\nu)$ is a Kazhdan $G$-space if and only if
$(X,\mu)$ is a Kazhdan $G$-space.
\end{thm*}

\begin{cor*}[\ref{zimmer}] Let $(Y,\nu)$ be an ergodic $G$-space where
$\nu$ is a finite $G$-invariant measure. Then the $G$-space $(Y,\nu)$
is Kazhdan if and only if $G$ is a Kazhdan group.
\end{cor*}

This extends Proposition 2.4 of Zimmer in \cite{zi:cohom} in one direction since,
in Zimmer's paper, the fact that $G$ is a Kazhdan group when the action has property
$T$ is only established under the additional assumption that the action is weakly
mixing. On the other hand, as it is well known, the assumption that $\nu$ is
finite and
$G$-invariant is crucial. For instance the action of a locally group $H$ on the
homogeneous space $H/K$ has property $T$ when the closed subgroup $K$ is Kazhdan,
although $H$ is not always Kazhdan in case $K$ is not of finite covolume.

As an immediate outcome of our approach, we get a characterization of property
$T$ for $r$-discrete measured groupoids by an appropriate boundedness condition
relative to  conditionally negative definite functions (Theorem
\ref{cond-neg-T}). It is used to deduce the fact, due to Adams and Spatzier
\cite{as} for equivalence relations, that an ergodic $r$-discrete measured groupoid in not
treeable in non trivial cases (see Corollary \ref{notreeable}).
Another consequence of this last characterization is given in the following theorem
(see Definition \ref{Haag-prop} for the Haagerup property).

\begin{thm*}[\ref{haag}] Every homomorphism from an ergodic Kazhdan $r$-discrete measured
groupoid to a locally compact group $H$ having the Haagerup property (i.e.,
a-$T$-menable in the sense of Gromov) is cohomologous to a homomorphism taking values in a compact subgroup
of $H$.
\end{thm*}

For an ergodic action of a Kazhdan locally compact group $G$ on $(X,\mu)$,
such that $\mu$ is finite and $G$-invariant, the corresponding theorem is
due to Jolissaint \cite[Theorem 3.2]{jo}. It extends previous results of \cite{as}.
\vspace{5mm}

This paper is organized as follows. In Section 2, we recall basic definitions on
bundles and groupoids. In Section 3, we define the first cohomology group
with respect to a representation and prove our main tool for the rest of the paper,
namely the characterization of coboundaries. Property $T$ is studied in Section 4,
and in particular the cohomological characterization is established.
Finally, Section 5 gives several examples of applications, and among them,
those described above.

To avoid further technicalities, we have limited our study to $r$-discrete
groupoids. The general case might be tackled using methods of \cite{fhm} and
\cite{ra:topo}.

In the whole paper, locally compact spaces are assumed to be second countable,
Borel spaces are standard, and Hilbert spaces are separable.
  
  {\it Acknowledgement.} After having written a preliminary version of this paper, I learned that
  some of the applications given in Section 5 are already contained in the unpublished thesis 
  (in Hebrew) of A. Nevo, in case of group actions and with a different approach. I am
grateful to him for several comments I took into account in this present version.

\section{Preliminaries}
\subsection{Borel bundles}A {\it bundle over a Borel space} $X$ is a Borel space
space $Z$ equipped with a Borel surjection $p: Z \ra X$. The Borel spaces $Z(x) =
p^{-1}(x)$ are the fibres of the bundle. Given an other bundle
$p': Z' \ra X$ over $X$, the {\it fibred product} of $Z$ and $Z'$ over $X$ is the
bundle
$$Z*Z' = \{(z,z') \in Z\times Z' : p(z) = p'(z')\}$$ endowed with the natural
projection and Borel structure. In case of ambiguity, we shall write $Z\,
_{p}\!*_{p'} Z'$ instead of $Z*Z'$.

A {\it Borel system of measures for} $p$ is a family $\rho = \{\rho^x : x\in X\}$ of
measures  $\rho^x$ on $p^{-1}(x)$ such that for every non-negative Borel function $f$
on $Z$, the function $x \mapsto \rho^x(f)$ is Borel. If there exists a positive Borel
function $f$ on $Z$ such that $\rho(f)$ is identically one, we say that $\rho$ is
{\it proper}.

 A {\it section} of the bundle $p: Z\ra X$ is a Borel map $\xi : X \ra Z$ with
$p\circ\xi(x) = x$ for all $x$.

\begin{defn}\label{metri} Let $X$ be a Borel space. A {\it bundle of metric spaces
over} $X$ is a bundle $p: Z\ra X$ equipped with a family 
$(d_x)_{x\in X}$  of
metrics $d_x$ on $Z(x)$ satisfying the two following conditions:
\begin{itemize}
\item[(1)] the map $(z,z') \mapsto d_{p(z)}(z,z')$ defined on $Z*Z$ is Borel.
\item[(2)] there exists a sequence $(\xi_n)$ of sections such that for every
$x\in X$, the set $\{\xi_n(x) : n\in \N\}$ is dense in $Z(x)$.
\end{itemize}
\end{defn}

As a consequence of the above condition $(1)$, note that for every section $\xi$
of $Z$, the map $(x,v) \mapsto d_x(v,\xi(x))$ is Borel on $X*Z$.

Let now $\h = \{\h(x)\}_{x \in X}$ be a family of Hilbert spaces  indexed by a
Borel set
$X$ and denote by $p$ the projection from $X*\h = \{(x,v) :v \in \h(x)\}$ 
to $X$.

\begin{defn}[\cite{ra:virt}, p. 264] A {\it Hilbert bundle} on a Borel space $X$ is a
space
$X*\h$ as above, endowed with a Borel structure such that
\begin{itemize}
\item[(1)] a subset $E$ of $X$ is Borel if and only if $p^{-1}(E)$ is Borel;
\item[(2)] there exists a {\it fondamental sequence $(\xi_n)$ of sections} satisfying
the following conditions :
\begin{itemize}
\item[(a)] for every $n$, the map $(x,v) \rightarrow \langle \xi_n(x), v\rangle$ is
Borel on $X*\h$;
\item[(b)] for every $m,n$, the map $x\rightarrow \langle \xi_m(x), \xi_n(x)\rangle$ is Borel;
\item[(c)] the functions  $(x,v) \rightarrow \langle \xi_n(x) , v\rangle$ separate
the points of $X*\h$.
\end{itemize}
\end{itemize}
\end{defn}
 
In Sections 2 and 3 we consider equally well either real or complex Hilbert bundles.

For Hilbert bundles we shall use the notation $X*\h$ or simply  $\h$. 
Note that every Hilbert bundle is a bundle of metric spaces, when the fibres
$\h(x)$ are endowed with the Hilbert metric. Let us also recall that when
all the fibres $\h(x)$ have the same dimension, the bundle is isomorphic to
a trivial Hilbert bundle (see \cite{ra:virt} for instance).

Given a Hilbert bundle $\h$ over $X$ and a class $C$ of measures on $X$, we denote by
$S((X,C), \h)$  the vector space of equivalence classes of Borel sections of the
Hilbert bundle $\h$,  where two sections are equivalent if they are equal
$C$-almost everywhere. We define a topology on  $S((X,C), \h)$ in the following way.
Having chosen a finite measure $\nu$ in $C$, we put
$$\rho(\xi, \eta) = \int_G \frac{\|\xi(\gamma) - \eta(\gamma)\|}
{1+\|\xi(\gamma) - \eta(\gamma)\|} d\nu(\gamma).$$
Obviously, $\rho$ is a metric on $S((X,C), \h)$.
Exactly as in \cite[Proposition 6]{mo:extIII}, one easily gets the following
result.

\begin{prop}\label{topolo} For a sequence $(\xi_n)$ in $S((X,C), \h)$ the
following are equivalent:
\begin{itemize}
\item[(i)] $\xi_n \rightarrow \xi$ for the metric $\rho$;
\item[(ii)] $\xi_n \rightarrow \xi$ in $\nu$-measure, that is, 
one has $\displaystyle \lim_{n\rightarrow \infty} \nu\Big(\frac{\|\xi_n -
\xi\|} {1+\|\xi_n - \xi\|} \geq \varepsilon\Big) =0$ for every
$\varepsilon >0$.
\item[(iii)] every subsequence of $(\xi_n)$ has a subsequence converging in norm
$C$-almost everywhere to $\xi$.
\end{itemize}
\end{prop}

It follows that the topology defined on $S((X,C), \h)$ only depends  on the measure
class $C$ of
$\nu$.

\begin{prop}\label{metr-comp} $S((X,C), \h)$ is a separable complete metrizable
topological vector space.
\end{prop}

\begin{proof} That $S((X,C), \h)$ is a topological vector space is easily checked. In
fact
$\xi \mapsto \displaystyle \int_G \frac{\|\xi\|}{1+\|\xi\|} d\nu$ is a $(F)$-norm
as defined in \cite[page 163]{ko}. Let us prove that
$\D\Big(S((X,C), \h), \rho\Big)$ is a complete metric space.  Let
$(\xi_n)$ be a Cauchy sequence. Taking if necessary a subsequence, we can assume
that
$\rho(\xi_{n+1}, \xi_n) \leq 1/2^n$ for every integer $n$. It follows that
$$\sum_n \frac{\|\xi_{n+1} - \xi_n\|}{1+ \|\xi_{n+1} - \xi_n\|} < +\infty
\quad C-\hbox{a.e.}$$
Therefore the series $\sum  \|\xi_{n+1} - \xi_n\|$ converges $C$-a.e,
and $(\xi_n)$ converges $C$-a.e. to an element $\xi$ of $S((X,C), \h)$.
Using Lebesgue's dominated convergence theorem, we see that
$\lim_{n\rightarrow +\infty} \rho(\xi_n, \xi) = 0$.
\end{proof}

Of course, $S((X,C), \h)$ is not a locally convex topological vector space in
general.

\subsection{Measured groupoids}
Our references for measured groupoids are \cite{ra:virt}, \cite{mu} and
\cite{dr:amenable}. Let us introduce first some notations. Given a groupoid
$G$, $G^{(0)}$ will denote its unit space and $G^{(2)}$ the set of composable pairs.
The range and source maps from $G$ to  $G^{(0)}$ will be denoted respectively
by $r$ and $s$. The corresponding fibres are denoted respectively $G^x =
r^{-1}(x)$ and $G_x = s^{-1}(x)$. Given subsets $A, B$ of $G^{(0)}$, we define $G^A =
r^{-1}(A)$,
$G_B = s^{-1}(B)$ and $G_{B}^A = G^A \cap G_B$. For $x\in G^{(0)}$, the isotropy
group $G^{x}_x$ will be denoted by $G(x)$. The
{\it reduction} of $G$ to
$A$ is the groupoid $G|_A = G_{A}^A$. Two elements $x,y$ of $G^{(0)}$ are said to be
equivalent if
$G^{x}_y \not = \emptyset$. If $A \subset G^{(0)}$, its {\it saturation} $[A]$ is the
set $s(r^{-1}(A))$ of
all elements in $G^{(0)}$ that are equivalent to some element of $A$. When $A = [A]$,
we say that $A$ is {\it invariant}.

From now on, $G$ will always be a {\it Borel groupoid}, which means that $G$ has a
Borel structure, for which the relevant operations are Borel maps. A {\it (left)
$G$-space} is a Borel space $X$ endowed with a Borel surjection $r : X \ra G^{(0)}$
and a Borel map $(\ga,x) \mapsto \ga x$ from the space $G \,_{s}\!*_r X$ into $X$
such that 
\begin{itemize}
\item[(i)] $r(\ga x) = r(\ga)$ for $(\ga, x) \in G\, _{s} \!*_r X$ , $r(x) x = x$ for
$x \in X$;
\item[(ii)] if $(\ga_1,x) \in G\,_{s}\!*_r X$, $(\ga_2,\ga_1) \in G^{(2)}$, then
$(\ga_2 \ga_1) x = \ga_2(\ga_1 x)$.
\end{itemize}
The corresponding {\it semi-direct groupoid} will be denoted $X\croi\, G$. 
Recall that $X\croi\, G = \{(x,\ga) : r(x) = r(\ga)\}$; the range and source maps
are respectively $r(x,\ga) = x$, and $s(x,\ga) = \ga^{-1} x$, and the product
is given by the formula $(x,\ga_1)(\ga_{1}^{-1}x ,\ga_2) = (x,\ga_1\ga_2)$.
One can define similarly right $G$-spaces. We sometimes view the unit space of a
Borel groupoid $G$ as a $G$-space, the action being given by $(\ga, s(\ga)) \mapsto
r(\ga)$.

A {\it homomorphism $\varphi: G\ra H$ between Borel groupoids} is a Borel map
that preserves the operations in an obvious way. Later, we shall sometimes call
$\varphi$ a {\it strict homomorphism} to avoid confusion with weaker notions of
homomorphisms.

Let $Z$ and $X$ be Borel $G$-spaces and let $p: Z \ra X$ be a $G$-equivariant
Borel surjection.
We say that a Borel system of measures $\rho$ for $p$ is {\it invariant} if $\gamma
\rho^x =\rho^{\ga x}$ whenever 
$s(\ga) = r\circ p(x)$, where 
$$\int f d\ga \rho^x = \int f(\ga z) d\rho^x(z).$$

\begin{defn}\label{haar} A {\it  Haar system} on $G$ is an invariant proper
Borel system $\lambda$ of measures for $r : G \ra G^{(0)}$.
\end{defn}

 Given a Borel groupoid $(G,\lambda)$ with Haar system, a measure
$\mu$ on $G^{(0)}$ is called {\it quasi-invariant with respect to} $(G,\lambda)$
if the measure $$\mu\circ\lambda : f \mapsto \int f(\ga)d\lambda^x(\ga) d\mu(x)$$
is equivalent to its image under $\ga \mapsto \ga^{-1}$.

\begin{defn}\label{sym-inv} Let $G$ be a Borel groupoid and $C$ a measure class on
$G$. 
\begin{itemize}
\item[(1)] We say that $C$ is {\it symmetric} if any measure $\nu \in C$ is
equivalent to its image $\nu^{-1}$ under inversion.
\item[(2)] We say that $C$ is {\it invariant} if there exists a probability measure
$\nu \in C$ whose disintegration with respect to $r : G \ra G^{(0)}$ is such that
$\ga \nu^{s(\ga)}$ is equivalent to $\nu^{r(\ga)}$ almost everywhere.
\end{itemize}
\end{defn}

\begin{defn}\label{meas-group} A {\it measured groupoid} is a pair $(G,C)$
such that $C$ is a symmetric invariant measure class on the Borel groupoid $G$.
In case $G$ is the semi-direct groupoid $X\croi\, H$ defined by the action of a Borel
groupoid $H$ with Haar system $\lambda$, together with the class of a quasi-invariant
measure on $X$, we also say that $X$ is a {\it measured $H$-space}. Note
that the Haar system for $X\croi\, H$ is then, implicitely, $\{\delta_x \times
\lambda^x : x\in X\}$, where $\delta_x$ is the Dirac measure at $x$.
\end{defn}

We shall denote by $r(C)$ the measure class of the image $r_*(\nu)$
of any probability measure $\nu\in C$.

\begin{defn} Let $(G,C)$ be a measured groupoid, and let $U\subset G^{(0)}$
be a subset that is conull with respect to $r(C)$. If $C|_U$ denotes the restriction
of the  measures in $C$ to $G|_U$, then $(G|_U, C|_U)$ is a measured groupoid, that
is called the {\it inessential reduction} of $G$ to $U$.
\end{defn}

\begin{defn} \label{ergo} A measured groupoid $(G,C)$ is called {\it ergodic}
in case the measure class $r(C)$ is ergodic in the sense that every invariant Borel
subset of $G^{(0)}$ is either null or conull for $r(C)$.
\end{defn}

An ergodic measured groupoid is a virtual group in the terminology of Ramsay
\cite{ra:virt}, following Mackey. 

Note that if $(G,\lambda)$ is a Borel group
with Haar system and $\mu$ is a measure on $G^{(0)}$ that is quasi-invariant with
respect to
$(G,\lambda)$, then $(G,C)$, where $C$ is the measure class of $\mu\circ\lambda$,
is a measured groupoid. Conversely, Hahn showed in \cite{ha1} that every measured
groupoid $(G,C)$ has an inessential reduction $(G|_U, C|_U)$, such that $C|_U$
is the class of a measure that can be decomposed in an essentially unique way as
$\mu\circ\lambda$, where
$\lambda$ is a Haar system for $G|_U$ and $\mu$ is a quasi-invariant measure on $U$
with respect to $(G|_U, \lambda)$. Therefore, the two point of views that are found
in the litterature, namely to define a measured groupoid as $(G,C)$ with $C$
symmetric and invariant, or as $(G,\lambda,\mu)$ where $\mu$ is quasi-invariant with
respect to the Haar system $\lambda$, are essentially the same. In the following, we
shall use the point of view which is more convenient, according to the
context. When the Haar system $\lambda$ is implicit, we shall sometimes denote by $(G, [\mu])$,
  or even by $(G,\mu)$,  the measured groupoid $G$
equipped with the class of the measure $\mu\circ\lambda$,   $[\mu]$
being the class of the measure $\mu$.

Given $\mu\circ\lambda \in C$, the class of the measure 
$$f\mapsto \int f(\ga_{1}^{-1},\ga_2) d\lambda^x(\ga_1)d\lambda^x(\ga_2) d\mu(x),$$
where $f$ is non-negative Borel on $G^{(2)} = \{(\ga_1,\ga_2) : s(\ga_1) =
r(\ga_2)\}$, will be denoted by $C^{(2)}$.

 Since we shall only be  interested in properties up to null sets, we
shall freely replace a measured groupoid by any of its inessential reductions, if
necessary. We shall also use repeatedly the following result.

\begin{lem}[Lemma 5.2, \cite{ra:virt}] \label{conull} Let $(G,C)$ be a measured
groupoid, let $H$ be a subset of $G$ that is closed under multiplication, and assume
that $H$ contains a conull set for $C$. Then there is a conull subset $U$ of
$G^{(0)}$ with respect to $r(C)$, such that $G|_U$ is contained in $H$.
\end{lem}

An important technical result, also due to Ramsay, says that we can always assume $G$
to be a
$\sigma$-compact topological groupoid. More precisely, we have:

\begin{thm}[Theorem 2.5, \cite{ra:topo}]\label{compact} Let $(G,C)$ be a measured
groupoid. Then there is
 an inessential reduction $(G|_U, C|_U)$ having a $\sigma$-compact
metric topology for which it is a topological groupoid and which defines the Borel
structure.
\end{thm}

In this context, the following lemma will be very useful for us.

\begin{lem}[\cite{ra:topo}]\label{federer} Let $(G,C)$ be a $\sigma$-compact
measured groupoid, and let $E$ be a subset of $G^{(0)}$ of positive measure. There
exists a conull subset $E_1$ of $E$ and a Borel section $\theta : [E_1] \ra
r^{-1}(E_1)$ for $s$.
\end{lem}

\begin{proof} Let us sketch the proof (see \cite{ra:topo} for details).
One takes $E_1$ to be $\sigma$-compact. Then $r^{-1}(E_1)$ is $\sigma$-compact,
as well as $[E_1] = s(r^{-1}(E_1))$, and one uses a reformulation of the
Federer-Morse Borel selection lemma.
\end{proof}

For technical reasons, in many cases we shall have to consider groupoids with
countable fibres.

\begin{defn}\label{rdiscret} We say that a Borel groupoid $G$ is $r$-{\it discrete}
if each fibre $G^x$ is countable. 
\end{defn}
 
Then the counting measures on the fibres
$G^x$ form a Haar system. Implicitly, $G$ will always be equipped with this Haar
system. A {\it measured groupoid $(G,C)$ is called $r$-discrete} in case the Borel
groupoid underlying some inessential reduction is $r$-discrete.

 Let us recall  two basic examples of measured groupoids.

\begin{ex}[{\it Discrete measured equivalence relations} \cite{fm1}] It
is an equivalence relation $R$ on a Borel space $X$ which has countable equivalence
classes, a Borel graph $R \subset X\times X$ and a quasi-invariant measure $\mu$.
This groupoid is $r$-discrete. The choice of the Haar system is therefore given by
the counting measures on the equivalence classes.  The measure $\mu$
is quasi-invariant in this case if for every Borel set $A\subset X$, $\mu(A) = 0$
implies $\mu([A]) = 0$. The
range
$r$ and source
$s$ are respectively the first and second projection from $R$ to $X$, the
composition is given by 
$(x,y)(y,z) = (x,z)$, and the inverse of $(x,y)$ is $(y,x)^{-1}$.
\end{ex} 

\begin{ex}[{\it Measured $G$-spaces}] Let $G$ be a locally compact group acting to
the left on a Borel space $X$, so that the left action map $X\times G \ra X$, $(x,g)
\mapsto gx$ is Borel. The semi-direct groupoid $X\croi\, G$ has already been
described previously for a groupoid action. The operations defining the structure
are $r(x,g) = x$, $s(x,g) = g^{-1}x$, $(x,g)(g^{-1}x, h) = (x,gh)$, $(x,g)^{-1}
= (g^{-1}x, g^{-1})$. The canonical Haar system is $\{\delta_x\times \lambda : x\in
X\}$, where  $\lambda$ is the left Haar
measure of $G$. A  measure on $X$ is quasi-invariant under the $G$-action in the
usual sense if and only if it is
quasi-invariant for the groupoid with Haar system $X\croi\, G$. Similarly, the two
notions of ergodicity coincide. 
\end{ex}

\section{First cohomology group with respect to a representation}

\subsection{Representations of a measured groupoid}

Given a Hilbert bundle $\h$ over a Borel space $X$, we denote by
$\hbox{Iso}(G^{(0)}*\h)$ the groupoid formed by the triples $(x,V,y)$, where
$x,y\in X$ and $V$ is a Hilbert isomorphism from $\h(y)$ onto $\h(x)$,
the composition law being defined by $(x,V,y)(y,W,z) = (x, V\circ W, z)$. We endow 
$\hbox{Iso}(G^{(0)}*\h)$ with the weakest Borel structure such that
$(x,V,y) \mapsto \l V\xi_n(y), \xi_m(x)\r$ is Borel for every $n,m$, where
$(\xi_n)$ is a fondamental sequence.

\begin{defn} A {\it representation of a Borel groupoid} $G$ is a pair
$(G^{(0)}*\h, L)$ where $G^{(0)}*\h$ is a Hilbert bundle over $G^{(0)}$,
and $L : G \rightarrow \hbox{Iso}(G^{(0)}*\h)$ is a Borel homomorphism that 
preserves the unit space. For $\gamma\in G$, we have $L(\gamma)
= (r(\gamma), \hat{L}(\gamma), s(\gamma)$, where $\hat{L}(\gamma)$ is an isometry
from
$\h (s(\gamma))$ onto  $\h(r(\gamma))$. To lighten the notations, we shall
identify $L(\gamma)$ and $ \hat{L}(\gamma)$. In particular we have
\begin{equation}\label{repre}
\forall (\gamma_1, \gamma_2) \in G^{(2)},\, L(\gamma_1\gamma_2) = L(\gamma_1)L(\gamma_2);\quad
\forall \gamma\in G,\, L(\gamma^{-1}) = L(\gamma)^{-1}.
\end{equation} 
\end{defn}

If the bundle $\h$ is trivial, with fibre $\ck$, then $L$ will be viewed as a Borel
homomorphism from $G$ into the unitary group $\mathcal{U}(\ck)$ equipped with the
Borel structure determined by the weak operator topology. \footnote{For simplicity
of language, we often use the term unitary instead of orthogonal, even in case of
real Hilbert spaces.}

\begin{defn} A {\it representation of a
measured groupoid} $(G,C)$ is a pair
$(G^{(0)}*\h, L)$ where $G^{(0)}*\h$ is a Hilbert bundle over
$G^{(0)}$, and $L : G \rightarrow \hbox{Iso}(G^{(0)}*\h)$ is a
Borel map such that there exists an inessential reduction
$G|_U$ for which the restriction of $L$ to $G|_U$ is a 
representation of the Borel groupoid $G|_U$.
\end{defn}

If we want to insist that (\ref{repre}) holds everywhere, we shall sometimes speak of
{\it strict representation}.

\begin{defn} Let $(G^{(0)}*\h_i, L_i)$, $i=1,2$, be two representations (in both real
or complex Hilbert bundles) of the measured groupoid
$(G,C)$. Let $U_i\subset G^{(0)}$ be a conull set such that the restriction of $L_i$
to $G|_{U_i}$ is a  strict representation. Then $L_1$ and $L_2$ are called
{\it equivalent} if there exists a conull subset $U \subset U_1 \cap U_2$ and a Borel
map 
$V : ( G^{(0)}*\h_1)|_{U} \rightarrow ( G^{(0)}*\h_2)|_{U}$
which is a fibre preserving isomorphism, and satisfies
\begin{equation}\label{equiv-rep}
V(r(\gamma)) L_1(\gamma) = L_2(\gamma) V(s(\gamma))
\end{equation}
for $\gamma \in G|_U$.
\end{defn}

By replacing the measured groupoid $(G,C)$ by an appropriate inessential reduction,
we shall always assume that (\ref{repre}) holds everywhere for representations of
measured groupoids.

\subsection{First cohomology group with respect to a representation} 

If $p : Y \rightarrow X$ is a Borel map between Borel spaces $Y$ and $X$,
and if $\h$ is a Hilbert bundle over $X$, 
the pull-back bundle $\{(y,v) : v \in \h(p(y))\}$ will be denoted by $Y*_p \h$
(or $p^* \h$).

\begin{defn}\label{cocy}Let $(G^{(0)}*\h, L)$ be a representation of the
measured groupoid $(G,C)$. Let $b : G \rightarrow G*_r \h$ be a Borel section.
\begin{itemize}
\item[(1)] $b$ is called a {\it strict $L$-cocycle} if for every $(\gamma_1,
\gamma_2) \in G^{(2)}$, we have 
\begin{equation}\label{cocy-egal}
b(\gamma_1\gamma_2) = b(\gamma_1) + L(\gamma_1)
b(\gamma_2).
\end{equation}
\item[(2)]  $b$ is called a $L$-{\it  cocycle} if there exists an inessential
reduction $G|_U$ such that (\ref{cocy-egal}) holds for every $(\gamma_1, \gamma_2)
\in G^{(2)} \cap (G|_U \times G|_U)$.
\item[(3)]  $b$ is called an {\it almost everywhere {\rm(a.e. for short)}
$L$-cocycle} if
(\ref{cocy-egal}) holds almost everywhere with respect to $C^{(2)}$.
\end{itemize}
\end{defn}

Let us comment on the relations between these  notions of $L$-cocycle. 
First, a $L$-cocycle is an a.e $L$-cocycle. Second, we have the following lemma,
whose proof is a straighforward adaptation of the proof of \cite[Theorem
5.1]{ra:virt}.

\begin{lem}\label{ae-cocy} Let $b$ be an a.e. $L$-cocycle. Then there exists
a $L$-cocycle $\tilde{b}$ such that $\tilde{b} = b$ a.e. with respect to $C$.
\end{lem}

Exactly as in \cite[Theorem 3.2]{ra:topo}, one can show the existence of an
inessential reduction $(G|_U,C|_U)$ such that for every $L$-cocycle $b$ there
exists a strict $L$-cocycle on $G|_U$ which agrees with $b$ almost everywhere
on $G|_U$. The
interesting fact is that the inessential reduction does not depend on $b$.
It is based on Theorem \ref{compact} and Lemma \ref{federer}.

Let us denote by $Z^1((G,C),L)$ the set of equivalence classes of a.e. $L$-cocycles
where two cocycles are equivalent if they are equal $C$-almost everywhere.

\begin{defn} We say that a  $L$-cocycle is a {\it $L$-coboundary} if there
exists a Borel section  $\xi$ of $G^{(0)}*\h$ such that
$b(\gamma) = \xi\circ r(\gamma) - L(\gamma) \xi\circ s(\gamma)$ for $C$-almost every $\gamma$.
We say that two  $L$-cocycles $b$ and $b'$ are {\it cohomologous} if $b-b'$
is a coboundary.
\end{defn}

Given a section $\xi$ of $\h$, the coboundary $\ga \mapsto \xi\circ r(\gamma) -
L(\gamma) \xi\circ s(\gamma)$ will often be denoted by $c_L(\xi)$.

We denote by $B^1((G,C),L)$ the set of equivalence classes of 
$L$-coboun\-daries, and  by $H^1((G,C),L) = Z^1((G,C),L)/B^1((G,C),L)$ the quotient
group, that we call the {\it first cohomology group of $(G,C)$ with respect to the
representation} $L$.

Observe  that replacing $(G,C)$ by any of its inessential reduction does not
change $H^1((G,C),L)$.

\begin{rem}\label{affine} Let $L$ be a representation of $G$ in $\h$, and $b$ a strict $L$-cocycle. For
$\ga \in G$, let us define $\alpha(\ga): \h(s(\ga)) \ra \h(r(\ga))$ by $\alpha(\ga)v = L(\ga)v + b(\ga)$.
Obviously, $\alpha(\ga)$ is an affine isometry, and $\alpha(\ga_1 \ga_2) = \alpha(\ga_1)\alpha(\ga_2)$
for $(\ga_1,\ga_2) \in G^{(2)}$, thanks to the cocycle property of $b$. A section $\xi$ of $\h$ is {\it
invariant} by $\alpha$, that is $\alpha(\ga)\xi\circ s(\ga) = \xi\circ r(\ga)$ for every $\ga \in G$, if and
only if
$b = c_L(\xi)$.
\end{rem}

\begin{rem}\label{isom} Up to isomorphism, $H^1((G,C),L)$ only depends on the
equivalence class of $L$. More precisely, if $L_1,L_2$ are two equivalent
representations such that (\ref{equiv-rep}) holds, given a $L_1$-cocycle $b_1$, let
us set $b_2(\ga) =V\circ r(\ga) b_1(\ga)$ for $\ga \in G$. Then $b_2$ is a
$L_2$-cocycle and $b_1 \mapsto b_2$ induces an isomorphism from $H^1((G,C), L_1)$
onto $H^1((G,C),L_2)$.
\end{rem}

\begin{prop} $Z^1((G,C),L)$ is a closed vector subspace of $S((G,C), \h)$
equipped with the topology defined in Proposition \ref{topolo}.
\end{prop}

\begin{proof} Let $(b_n)$ be a sequence in $Z^1((G,C),L)$, converging to $\varphi \in S((G,C), L)$.
By Lemma \ref{ae-cocy}, we may assume the existence of an inessential reduction
$G|_U$ such that each $b_n$ restricted to $G|_U$ is a strict cocycle.  Observe that
the set
$\{ \gamma \in G|_U : (b_n(\gamma)) \,\,\hbox{converges}\}$ is stable by composition and
inverse. Moreover this set is conull in $G|_U$. By Lemma \ref{conull}, this set contains
an inessential reduction $G|_V$ of $G|_U$. It follows, since the $L$-cocycle
identities (\ref{cocy-egal}) remain fulfilled by passing to the limit, that
$\varphi$ restricted to
$G|_V$ is a strict $L$-cocycle, and therefore $\varphi \in Z^1((G,C),L)$.
\end{proof}

In general $B^1((G,C),L)$ is not closed and, as we shall see, the fact for
$B^1((G,C),L)$ to be closed has important consequences.

\subsection{Homomorphisms and functoriality} Cohomology of measured grou\-poids with
coefficients in a Borel group is a fundamental subject in ergodic theory,
which has been widely studied (see \cite{we}, \cite{mo:extIII}, \cite{fm1},
\cite{sch1} for instance). We just need here to introduce some terminology.

\begin{defn}\label{hom} Let $(G,C)$ be a measured groupoid and $H$ a Borel
groupoid. 
\begin{itemize}
\item[(i)] A {\it homomorphism} or {\it cocycle} from $(G,C)$ into $H$
is a Borel map $\varphi : G\ra H$ such that there exists an inessential reduction
$(G|_U, C|_U)$ on which the restriction of $\varphi$ is a (strict) homomorphism 
of Borel groupoids.
\item[(ii)] Two homomorphisms $\varphi, \psi$ are {\it cohomologous} if there exists
an inessential reduction $(G|_U, C|_U)$ on which they are both strict homomorphisms,
and a Borel map $\theta : U \ra H$ such that $\theta(r(\ga))\varphi(\ga)$ and
$\psi(\ga)\theta(s(\ga))$ are defined and equal for $\ga \in G|_U$.
\end{itemize}
\end{defn}

The notion of homomorphism between measured groupoids is much more
intricate. It was clarified by Ramsay in
\cite{ra:virt}. We shall not enter here into this subject, and shall limit
ourself to describe the situation in the case of $r$-discrete measured groupoids.

Given a homomorphism $\varphi : G_1 \ra G_2$ between two Borel groupoids,
we shall often denote by $\tilde{\varphi} : G_1^{(0)} \ra G_2^{(0)}$ its restriction
to the unit spaces.
 
\begin{defn}\label{meas-hom}  Let $(G_1,C_1)$ and
$(G_2,C_2)$ be two $r$-discrete measured grou\-poids.
\begin{itemize}
\item[(1)] A {\it strict homomorphism of $r$-discrete measured groupoids} $\varphi : 
(G_1,C_1) \ra (G_2,C_2)$  is a strict Borel homomorphism such that
$\tilde{\varphi}^{-1}(E)$ is null with respect to $r(C_1)$ for every 
$r(C_2)$-null Borel subset $E$ of $G_2^{(0)}$.
\item[(2)] A {\it homomorphism of $r$-discrete measured groupoids} $\varphi : 
(G_1,C_1) \ra (G_2,C_2)$ is a Borel map, for which there exists an
inessential reduction $(G_1|_U, C_1|_U)$ such that the restriction of $\varphi$ to
$(G_1|_U, C_1|_U)$ is a strict homomorphism of $r$-discrete measured groupoids.
\end{itemize}
\end{defn}

\begin{defn}  Let $(G_1,C_1)$ and $(G_2,C_2)$ be two $r$-discrete
measured grou\-poids.
\begin{itemize}
\item[(1)] A {\it strict similarity between two strict homomorphisms} $\varphi$, $\psi$ from
$(G_1,C_1)$ to $(G_2,C_2)$ is a Borel map $\theta : G_1^{(0)} \ra G_2$, such that
$\theta(r(\ga))\varphi(\ga)$ and $\psi(\ga)\theta(s(\ga))$ are defined and equal for $\ga \in G_1$. 
\item[(2)] For two homomorphisms $\varphi$, $\psi$, we say that $\varphi$ is {\it similar} 
(or {\it cohomologous}) to $\psi$
if there is an inessential reduction $(G_1|_U, C_1|_U)$ such that the restrictions of $\varphi$ and $\psi$ to
$(G_1|_U, C_1|_U)$ are strict and strictly similar. In this case, we write $\varphi \sim \psi$.
\end{itemize}
\end{defn}

In the langage of \cite{fhm}, a similarity between homomorphisms is called an
equivalence.

\begin{rem} For general measured groupoids, one only asks for a strict
homomorphim $\varphi$ to be such that $\tilde{\varphi}^{-1}(E)$ is null with 
respect to $r(C_1)$ for every saturated
$r(C_2)$-null Borel subset $E$ of $G_2^{(0)}$. Thus a $r(C_2)$-null set $E$ in
$G_2^{0)}$ is not always null for the image of $r(C_1)$ by $\tilde{\varphi}$, but it
is the case whenever $E$ is ``small'' in the sense that its saturation $[E]$ is
still $r(C_2)$-null.

One of the main motivations for this weaker notion, was the fundamental observation
made by Mackey that, if $H$ is a closed subgroup of a locally compact group $G$,
then the groupoid $G/H\croi\, G$ is strictly similar to the group
$G$ (see \cite[Definition 6.4, Theorem 6.20]{ra:virt} for the precise 
definition and the proof).
The principal difficulty when working with this weaker definition is the
composition of homomorphisms (see \cite{ra:virt}). 

For $r$-discrete groupoids
this weaker definition of homomorphism coincide with the notion introduced in
Definition \ref{meas-hom} because in this case a subset $E$ of the unit space is
null if and only if its saturation is null. There is of course no difficulty
to compose homomorphisms in this case.
\end{rem}

Now, let $\varphi : (G_1,C_1) \ra (G_2,C_2)$ be a homomorphism of $r$-discrete
measured groupoids, and consider a $L_2$-representation of $(G_2,C_2)$ in a Hilbert
bundle $\h_2$. By passing if necessary to inessential reductions, we can assume that
$\varphi$ and $L_2$ are strict. Let us denote by $\h_1 =\tilde{\varphi}^* \h_2$ the
pull-back of $\h_2$, that is
$$G_1^{(0)} * \h_1 = \{(x,v) : x\in G_1^{(0)}, v \in \h_2(\tilde{\varphi}(x))\}.$$
Then $L_1 = L_2\circ\varphi$ is a representation of $G_1$ in $\h_1$. If $b_2$ is a $L_2$-cocycle
then $b_2\circ\varphi$ is a $L_1$-cocycle. In this way, we define without ambiguity
a map from $Z^1((G_2,C_2),L_2)$ into $Z^1((G_1,C_1),L_1)$, that sends $B^1((G_2,C_2),L_2)$
into $B^1((G_1,C_1),L_1)$. Therefore we get a map
$$\varphi^* : H^1((G_2,C_2),L_2) \ra H^1((G_1,C_1),L_1).$$

We easily check that given another homomorphism $\psi : (G_0, C_0) \ra (G_1,C_1)$
of $r$-discrete measured groupoids, we have $(\varphi\circ\psi)^* = \psi^*\circ\varphi^*$. 

Let us compare now the maps $\varphi^*$ and $(\varphi')^*$ for two similar homomorphisms. We can assume that
$\varphi$ and $\varphi'$ are strictly similar strict homomorphisms of measured groupoids. So, there is
a Borel map $\theta : G_1^{(0)} \ra G_2$ such that 
$\theta(r(\ga))\varphi(\ga)$ and $\varphi'(\ga)\theta(s(\ga))$ are defined and equal for $\ga \in G_1$.
From $\varphi'$, we define the representation $L_1'$ of $G_1$ in $\h_1'$ as we did for $\varphi$. Note that
the representations $L_1$ and $L_1'$ are equivalent: $L_2(\theta(x))$ is an isomorphism from
$\h_1(x) = \h_2(\tilde{\varphi}(x))$ onto $\h'_{1}(x) = \h_2(\tilde{\varphi}'(x))$ for every $x\in G_1^{(0)}$ and
we have
$$L_2(\theta\circ r(\ga)) L_1(\ga) = L'_{1}(\ga) L_2(\theta\circ s(\ga)),$$
for $\ga \in G_1$. Let us denote by $\Theta$ the isomorphism from $H^1((G_1,C_1), L_1)$ onto
$H^1((G_1,C_1), L_1')$ defined by the map $b \mapsto b'$ with $b'(\ga) = L_2\circ\theta(r(\ga)) b(\ga)$
(see \ref{isom}).

\begin{lem}\label{comp-sim} We have $(\varphi')^* = \Theta\circ \varphi^*$.
\end{lem}

\begin{proof} This results from a straighforward computation. Let $b$ be a $L_2$-cocycle. For $\ga \in
(G_1)^{x}_y$, we have, using the cocycle identity,
\begin{align*}
b(\varphi'(\ga)) &= b(\theta(x)) + L_2(\theta(x))\Big( b(\varphi(\ga)) -
L_2(\varphi(\ga))L_2(\theta(y))^{-1} b(\theta(y))\Big)\\
&= L_2(\theta(x))b(\varphi(\ga)) + b(\theta(x)) - L'_{1}(\ga) b(\theta(y)),
\end{align*}
and therefore $\ga \mapsto b(\varphi'(\ga))$ and $\ga\mapsto L_2\circ\theta(r(\ga))b(\varphi(\ga))$
are cohomologous.
\end{proof}

\begin{defn}\label{simil} A homomorphism $\varphi : (G_1,C_1) \ra (G_2,C_2)$ of $r$-discrete measured groupoids
is called a {\it similarity}, and we say that the groupoids are {\it similar},  if there exists a homomorphism
$\psi : (G_2,C_2) \ra (G_1,C_1)$ such that $\varphi\circ \psi \sim \iota_{G_2}$ and
$\psi\circ\varphi  \sim \iota_{G_1}$, where $\iota_{G_1}$ and $\iota_{G_2}$ are the identity
homomorphisms of $G_1$ and $G_2$ respectively.
\end{defn}

\begin{prop}\label{stable} Let $\varphi : (G_1,C_1) \ra (G_2,C_2)$ be a similarity between $r$-discrete measured
groupoids. Given a representation  $L_2$ of $(G_2,C_2)$ in $\h_2$, we denote by $L_1 = L_2\circ \varphi$
the corresponding representation of $(G_1,C_1)$ in $\tilde{\varphi}^* \h_2$. Then the
induced application
$$\varphi^* : H^1((G_2,C_2), L_2) \ra H^1((G_1,C_1), L_1)$$
is an isomorphism.
\end{prop}

\begin{proof} Let $\psi$ as in Definition \ref{simil}, and define the representation $L_2' = L_1\circ\psi$ of
$(G_2,C_2)$. By Lemma \ref{comp-sim}, there is an isomorphism 
$$\Theta : H^1((G_2,C_2), L_2) \ra H^1((G_2,C_2),L'_2)$$
 such that $\Theta = \Theta\circ \iota_{G_2}^* =
\psi^*\circ \varphi^*$. It follows that $\varphi^*$ is injective and $\psi^*$ is surjective. Similarly, we get
that  $\psi^*$ is injective and  $\varphi^*$ is surjective.
\end{proof}

\subsection{A necessary and sufficient condition for a $L$-cocycle to be a coboun\-dary}
We shall adapt to our context a proof due to Serre when the groupoid is a locally
compact group. A good reference for this
subject is \cite{hv}. We shall say that a metric space $(Y,d)$ satisfies the {\it
median inequality} if the two following conditions are fulfilled:
\begin{itemize}
\item[(i)] for every $x,y \in Y$ there exists a unique point $m$ (called the
{\it middle of the segment} $[x,y]$) such that
$d(x,m) = d(m,y) = \frac{1}{2} d(x,y)$. 
\item[(ii)] given three points $a,b,c$ of $Y$, and denoting by $m$ the middle
of $[b,c]$, we have the median inequality 
$$2d(a,m)^2 +\frac{1}{2} d(b,c)^2 \leq d(a,b)^2 + d(a,c)^2.$$
\end{itemize}
For instance, Hilbert spaces, trees and real trees satisfy the {\it median
inequality}; for more examples see \cite[page 38]{hv}.

\begin{defn} Let $G$ be a Borel groupoid. A {\it $G$-bundle of metric spaces over} $X$ is a bundle
$p: Z \ra X$ of metric spaces such that $Z$ and $X$ are Borel $G$-spaces and $p$ is $G$-equivariant,
and where $G$ acts fibrewise by isometries. An {\it invariant section} is a Borel section $\xi$ such that 
$\ga\xi(x) = \xi(\ga x)$ for $(\ga,x) \in G\, _s\! *_r X$.
\end{defn}

For instance, if $L$ is a representation of $G$ in $\h$, then $G^{(0)}*\h$ is a $G$-bundle of metric spaces over 
$G^{(0)}$.

Our first result gives a sufficient condition for a $G$-bundle to have an invariant
section.

\begin{lem}[see Proposition 9, page 38, \cite{hv}]\label{inv-sect} Let $(G,C)$ be a
measured \newline groupoid. Let $Z$ be a Borel $G$-space such that $p: Z \ra
G^{(0)}$ is a
$G$-bundle over $G^{(0)}$ of complete metric spaces $(Z(x), d_x)$ satisfying the median inequality. 
Suppose we are given a family $(B(x))_{x\in G^{(0)}}$ of bounded subsets $B(x)$
of $Z(x)$ and assume that
\begin{enumerate}
\item $\gamma B(s(\ga)) = B(r(\ga))$ for $C$-almost every $\ga$;
\item there exists a sequence $(\eta_n)_{n\in \n}$ of Borel sections $x\mapsto
\eta_n(x)\in B(x)$, such that for $r(C)$-almost every $x$, the
set $\{\eta_n(x) : n\in \N\}$ is dense in $B(x)$.
\end{enumerate}
Then there exists an inessential reduction $(G|_U, C|_U)$ and a  Borel section
$x\in U \mapsto \xi(x) \in Z(x)$ such that $\ga \xi\circ s(\gamma) =
\xi\circ r(\ga)$ for every $\ga \in G|_U$.
\end{lem}

\begin{proof} 
By restriction to an inessential reduction, we may assume that the conditions in
(i) and (ii) hold everywhere. For
every $x\in G^{(0)}$ we define $\rho_x : Z(x)\ra\R^+$ by
$$\rho_x(v) = \sup_{w\in B(x)}d_x(v,w).$$
Since 
$$
\rho_x(v) = \sup_{n\in \n}d_x(v,\eta_n(x))$$
we see that the map $(x,v) \in G^{(0)}*Z \mapsto \rho_x(v)$ is Borel.
Note also that the map $v\mapsto \rho_x(v)$ is continuous on $Z(x)$.

Now we set $r(x) := \inf_{v\in Z(x)} \rho_x(v)$. Using the separability assumption
(2) in Definition \ref{metri},
it is easily checked that
$r$ is a Borel map on $G^{(0)}$.

For every integer $n \geq 1$, let us define
$$D_n = \{ (x,v) \in G^{(0)}*Z : \rho_x(v) \leq r(x) + \frac{1}{n}\}.$$
Then $D_n$ is a Borel subset of $G^{(0)}*Z$.
By the von Neumann selection theorem \cite[Theorem A.9, p. 196]{zi:ergo}, and
replacing if necessary
$(G,C)$ by an inessential reduction, we can choose a Borel section $\xi_n$ of the
bundle
$Z$ such that $\xi_n(x) \in D_n(x)$ for every $x \in G^{(0)}$.
Using the median inequality, as in \cite[pages 37-38]{hv}, we show
that $\big(\xi_n(x)\big)$ is a Cauchy sequence in $Z(x)$, for every $x\in G^{(0)}$.
Denote by $\xi(x)$ its limit. Observe that $\rho_{x}(\xi(x)) = r(x)$,
and that $\xi(x)$ is the unique element in $Z(x)$ where the minimum
of the function $\rho_x$ is achieved (see \cite[Lemme du centre, p.37]{hv}.
It is the centre of the ball of smallest radius containing $B(x)$.
The main point for us is that the map $\xi$ is Borel on $G^{(0)}$. 

Now, since $\ga$ induces an isometry from $Z(s(\ga))$ onto $Z(r(\ga))$,
and since it carries $B(s(\ga))$ onto $B(r(\ga))$, it sends the ``centre''
$\xi\circ s(\ga)$ of $B(s(\gamma))$ onto the ``centre''
$\xi\circ r(\ga)$ of $B(r(\gamma))$.
\end{proof}

\begin{thm}\label{triv-cocy} Let $(G,C)$ be a measured
groupoid, 
$(G^{(0)}*\h, L)$ a representation of $(G,C)$, and $b$ a $L$-cocycle. Let us consider
the following conditions:
\begin{itemize}
\item[(i)] $b$  is a $L$-coboundary;
\item[(ii)] there exists a Borel subset $E$ of $G^{(0)}$, of
positive measure, such that $\sup\{ \|b(\gamma)\| :
\gamma \in G|_E\} < +\infty$;
\item[(iii)] there exists a Borel subset $E$ of $G^{(0)}$, of
positive measure, such that 
for every $x\in E$, we have $\sup\{ \|b(\gamma)\| :\gamma \in G_{E}^x \} < +\infty$.
\end{itemize}
Then we have $(i) \Rightarrow (ii) \Rightarrow (iii)$ and  when $(G,C)$ is
an ergodic $r$-discrete measured groupoid, the three conditions are equivalent.
\end{thm}

\begin{proof} (i) $\Rightarrow$ (ii). Assume that there exists a Borel section $\xi$ of $\h$ such
that $b(\gamma) = \xi\circ r(\gamma) - L(\gamma)\xi\circ s(\gamma)$ for $C$-almost
every $\ga$. Let
$(G|_U, C|_U)$ be an inssential reduction where this equality holds everywhere. Since
$U =
\cup_{n\geq 1} E_n$, where $E_n = \{ x\in U : \| \xi(x)\| \leq n\}$, there exists 
an integer $n_0$ such that the measure of $E := E_{n_0}$ is positive. Then for $\gamma \in G|_E$
we have $$\|b(\gamma)\| = \| \xi\circ r(\gamma) - L(\gamma)\xi\circ s(\gamma)\| \leq 2 n_0.$$

(ii) $\Rightarrow$ (iii) is obvious and the fact that (iii) $\Rightarrow$ (i)
when when $(G,C)$ is an ergodic $r$-discrete measured groupoid will be a consequence
of the two following lemmas.
\end{proof}

\begin{lem}\label{boun-cocy} Let $b$ be a $L$-cocycle of an ergodic groupoid $(G,C)$.
Assume the existence a Borel subset $E$ of $G^{(0)}$, of positive measure,
such that for every $x\in E$, we have $\sup\{ \|b(\gamma)\| : \gamma \in G_{E}^x \} < +\infty$.
Then $b$ is cohomologous to a $L$-cocycle $b'$ such that there exists 
a conull subset $U$ in $G^{(0)}$ with $\sup\{ \|b'(\gamma))\| : \gamma \in G_{U}^x \}
< +\infty$ for every $x\in U$.
\end{lem}

\begin{proof} We make several simplifications by passing to inessential
reductions. First, we may assume that $b$ is a strict $L$-cocycle and, thanks to
\ref{compact}, that
$G$ is a $\st$-compact topological groupoid. Moreover, since $(G,C)$
is ergodic, we may suppose that the bundle
$\h$ is trivial,  so that $L$ is a homomorphism from $G$ into the unitary
group of a separable Hilbert space $\ck$.
Finally, replacing if necessary $E$ by a conull subset, we can assume, by Lemma \ref{federer}, the existence of
a Borel section $\theta : [E] \ra r^{-1}(E)$ of $s$. Now we define $b_{\theta}: G \ra \ck$ as follows:
\begin{align*}
b_{\theta}(\gamma) &= L\big(\theta\circ r(\gamma)\big)^{-1} b\big(\theta\circ
r(\gamma)\gamma (\theta\circ s(\gamma))^{-1}\big) \quad \hbox{if}\quad
\gamma\in G|_{[E]};\\
b_{\theta}(\gamma) &= 0 \quad \hbox{if}\quad\gamma\notin G|_{[E]}.
\end{align*}
Observe that $\theta\circ r(\gamma)\gamma (\theta\circ s(\gamma))^{-1} \in
G^{r(\theta\circ r(\ga))}_E$ for every $\gamma \in G|_{[E]}$ so that, for $x \in
[E]$, the set
$$\{b_{\theta}(\ga) : \ga \in G^{x}_{[E]}\}$$ is  bounded in $\ck$ because $r(\theta(x)) \in E$.
Moreover, for
$\gamma\in G|_{[E]}$ we have, since $b$ is a $L$-cocycle,
\begin{align*}
b(\gamma) & = b\big((\theta\circ r(\gamma))^{-1}\theta\circ r(\gamma)\gamma\big)\\
&= b\big((\theta\circ r(\gamma))^{-1}\big) + L(\theta\circ r(\gamma))^{-1}
b\big(\theta\circ r(\gamma)\gamma\big),
\end{align*}
and 
$$
b\big(\theta\circ r(\gamma)\gamma\big)= b\big(\theta\circ
r(\gamma)\gamma (\theta\circ s(\gamma))^{-1}\big) + 
L\big(\theta\circ r(\gamma)\gamma (\theta\circ s(\gamma))^{-1}\big)
b(\theta\circ s(\gamma)),$$
so that
\begin{align*}
b(\gamma) &= b_{\theta}(\gamma) + b((\theta\circ r(\gamma))^{-1})
- L(\gamma) b((\theta\circ s(\gamma))^{-1})\\
&= b_{\theta}(\gamma) + \xi\circ r(\gamma) -L(\gamma)\xi\circ s(\gamma),
\end{align*}
if $\xi$ denotes the map $x\mapsto b(\theta(x)^{-1})$ on $[E]$.
Hence, $b$ and $b_\theta$ are cohomologous.
\end{proof}

\begin{lem}\label{boun-cocy2} Let $(G,C)$ be a $r$-discrete measured
groupoid and $(G^{(0)}*\h, L)$ a representation. Let $b$ be a $L$-cocycle
such that for every $x\in G^{(0)}$, the set $B(x) :=\{b(\ga) : \ga \in G^x\}$ is
bounded in $\h(x)$. Then $b$ is a coboundary.
\end{lem}

\begin{proof}  For $\gamma\in G$,
let  $\alpha(\gamma) : v\mapsto  L(\gamma)v + b(\gamma)$ be the corresponding
affine representation.
Then $\alpha$ defines on $\h$ a structure of $G$-bundle of metric spaces over $G^{(0)}$
which are complete and satisfy the median inequality. By the cocycle
property of $b$, we have $\alpha(\ga)\big(B(s(\ga))\big) = B(r(\gamma))$.
On the other hand, since $G$ is $r$-discrete, there exists a countable family
$(\sigma_n)$ of Borel sections 
$$\sigma_n : x\in G^{(0)} \mapsto \sigma_n(x) \in G^x$$ 
of $r$, such that $G = \cup
\sigma_n(G^{(0)})$ and so $B(x) = \cup_{n} B(\sigma_n(x))$.  Then, it follows from
Lemma \ref{inv-sect} that there is a Borel section $\xi$ of the Hilbert bundle $\h$
such that $\alpha(\ga) \xi\circ s(\ga) = \xi\circ r(\ga)$ almost everywhere, and
therefore
$b(\gamma) = \xi\circ r(\gamma) - L(\ga) \xi\circ s(\gamma)$
almost everywhere.
\end{proof}

\section{Property T}

From now on, unless specific indications are given, we implicitly consider representations
in complex Hilbert bundles.

\subsection{Definition of property T}
Let us denote by $\hbox{Rep}(G,C)$ the space of representations of $(G,C)$
(in complex Hilbert bundles),
where equivalent representations are identified. We first define a topology on
$\hbox{Rep}(G,C)$, which is similar to the topology defined by Fell
when $G$ is a group \cite{fe}.

Given a Hilbert bundle $\h$ on $G^{(0)}$, we shall denote by $S_1((G^{(0)},r(C)),\h)$
the space of all its Borel sections $\xi$ such that $\|\xi(x)\| = 1$ almost
everywhere, where, as always, two sections which agree almost everywhere
are identified. Elements of $S_1((G^{(0)},r(C)),\h)$ are called {\it unit sections}\footnote{
If the bundle $\h$ is trivial with fibre $\ck$, we shall often write 
$S_1((G^{(0)},r(C)),\ck)$ instead of $S_1((G^{(0)},r(C)),\h)$}.

Let $(G^{(0)}*\h, L) \in \hbox{Rep}(G,C)$. For
$\xi,\xi' \in S_1((G^{(0)},r(C)),\h)$, the corresponding {\it coefficient of the
representation} is the map $(\xi,\xi')_L$ defined on $G$ by
$$(\xi,\xi')_L(\ga) = \l \xi\circ r(\ga), L(\ga) \xi\circ s(\ga)\r.$$ 

Given a neighbourhood
$\mathcal{V}$ of $0$ in $L^\infty(G,C)$ for the weak*-topology, and
a finite set $\{\xi_1,\dots,\xi_n\}$ of elements of $S_1((G^{(0)},r(C)),\h)$,
we denote by $$V(L; \xi_1,\dots,\xi_n, \mathcal{V})$$ the set of representations
$(G^{(0)}*\h', L') \in \hbox{Rep}(G,C)$ such that there exists
$\eta_1,\dots,\eta_n \in S_1((G^{(0)},r(C)),\h')$ with
$$(\xi_i,\xi_j)_L - (\eta_i, \eta_j)_{L'} \in \mathcal{V}, \quad 1\leq i,j\leq n.$$ 
Then $\hbox{Rep}(G,C)$ is endowed with the topology such that, for every
representation $(G^{(0)}*\h, L)$, the sets $V(L;\dots)$ form a basis
of neighbourhoods of $(G^{(0)}*\h, L)$.

We are mainly interested in neighbourhoods of the {\it trivial representation} $L_0$
of $(G,C)$, and in this situation there are simpler basis of neighbourhoods. Recall
that $L_0$ is the representation in the trivial bundle over $G^{(0)}$
with fibre $\C$, such that $L_0(\ga) =1$ for $\ga \in G$. To define
a basis of neighbourhoods of $L_0$ it is enough to consider the constant
section equal to $1$ in $S_1((G^{(0)},r(C)),\C)$. Then, for $\mathcal{V}$ as above,
$V(L_0, \mathcal{V})$ will be the set of representations
$(G^{(0)}*\h', L') \in \hbox{Rep}(G,C)$ such that there exists
$\eta \in S_1((G^{(0)},r(C)),\h')$ with
$1 - (\eta, \eta)_{L'} \in \mathcal{V}$. These  sets $V(L_0, \mathcal{V})$
still form a basis of neighbourhoods of $L_0$.

There are other variants, that are even more convenient. Given  a representation
$(G^{(0)}*\h, L)$  of $(G,C)$ and a section 
$\xi$ of the bundle $G^{(0)}*\h$, we recall that we denote by $c_L(\xi)$
the coboundary $\ga \mapsto \xi\circ r(\ga) - L(\ga)\xi\circ s(\ga)$.
Let $W(L_0, \mathcal{V})$ be the set of representations
$(G^{(0)}*\h, L) \in \hbox{Rep}(G,C)$ such that there exists
$\xi \in S_1((G^{(0)},r(C)),\h)$ with 
$\|c_L(\xi)\| \in \mathcal{V}$. It follows immediately from the Cauchy-Schwarz
inequality that these sets $W(L_0, \mathcal{V})$ also form a
basis of neighbourhoods of $L_0$. This is classical for group representations, and
works as well in our context. One uses, for $\xi \in S_1((G^{(0)},r(C)),\h)$,
 the inequalities
$$|1- (\xi,\xi)_L| \leq \|c_L(\xi)\| \leq \sqrt{2} |1- (\xi,\xi)_L|^{1/2}.$$

\begin{lem}\label{equi-sect} Let $(G^{(0)}*\h, L)$ be a representation of $(G,C)$ and
let $\nu$ be a probability measure in the class $C$.
The following conditions are equivalent:
\begin{enumerate}
\item for every $\varepsilon >0$, there exists $\xi \in S_1((G^{(0)},r(C)),\h)$
such that 
$$\nu(\|c_L(\xi)\| \geq \varepsilon) \leq \varepsilon;$$
\item there exists a sequence of elements $\xi_n \in S_1((G^{(0)},r(C)),\h)$ such
that $\D\lim_{n\ra +\infty} c_L(\xi_n) = 0$ $C$-a.e.;
\item there exists a sequence of elements $\xi_n \in S_1((G^{(0)},r(C)),\h)$ such
that $\D\lim_{n\ra +\infty} c_L(\xi_n) = 0$ in $L^\infty(G)$ endowed with the
weak*-topology.
\end{enumerate}
\end{lem}

\begin{proof} (i) $\Rightarrow$ (ii) follows from the fact that every sequence
converging to $0$ in measure contains a subsequence converging to $0$ almost
everywhere. For the converse (ii) $\Rightarrow$ (i) one uses the fact that,
for a finite measure, the convergence a.e. implies the convergence in measure.
The equivalence between (ii) and (iii) is an immediate consequence
of the Lebesgue convergence theorems.
\end{proof}

\begin{defn}\label{invsect} Let $(G^{(0)}*\h, L)$ be a representation of $(G,C)$.
\begin{itemize}
\item[(1)] A section $\xi$ of the bundle $G^{(0)}*\h$ is said to
be {\it invariant with respect to} $L$ if $\xi\circ r(\ga) = L(\gamma)\xi\circ
s(\ga)$ a.e.
\item[(2)] If there is an invariant section $\xi \in S_1((G^{(0)},r(C)),\h)$
with respect to $L$,
we say that the representation $L$ {\it contains a unit invariant section}, or that
{\it the trivial representation $L_0$ of $(G,C)$ is contained into $L$}.
\item[(3)] We say that $L$ {\it almost contains unit invariant sections} if
the equivalent conditions of Lemma \ref{equi-sect} are fulfilled.
In this case, we also say that the trivial representation is {\it weakly
contained into} $L$.
\end{itemize}
\end{defn}

Note that the above property (2) means that the trivial representation
is a direct factor of $L$. When the groupoid is ergodic, this is equivalent
to saying that there is a non-zero invariant section for $L$. Property (3) means that
the trivial representation $L_0$ belongs to the closure of $L$ in $\hbox{Rep}(G,C)$.

\begin{defn}\label{prop-T} We say that the measured groupoid $(G,C)$
has {\it property} $T$,  or is a {\it Kazhdan groupoid}, if every representation
of $(G,C)$ which almost has unit invariant sections actually has unit invariant
sections.
\end{defn}

\begin{exs} When the groupoid is a locally compact group, we recognize the classical
definition of property
$T$ introduced by Kazhdan in his seminal paper \cite{ka}. For a group $G$, usually one says that the
trivial representation $L_0$ of
$G$ is weakly contained in a representation $(\pi, \h)$ if and only if there is a
sequence $(\xi_n)$ of unit vectors in $\h$ such that the sequence of functions $g
\mapsto \|\pi(g)\xi_n - \xi_n\|$ goes to zero uniformly on compact subsets of $G$
(see \cite[\S 18.3]{dix}. An equivalent formulation (\cite[\S 3.4]{dix} for instance)
is the existence of a sequence $(\xi_n)$ of unit vectors in $\h$  such that for every
$f\in L^1(G)$ we have
$$\lim_n \l \pi(f)\xi_n, \xi_n \r = \int f(g) dg,$$
that is, the sequence of coefficients $(\xi_n, \xi_n)_\pi$ goes to $1$ in the
weak*-topology of $L^\infty(G)$. This is the notion of weak containment defined in
\ref{invsect}. (3) above. In fact, our topology on $\hbox{Rep}(G)$ is the Fell
topology in case of groups.

For discrete measured equivalence relations, Definition \ref{prop-T} is the
definition introduced by Zimmer in \cite{zi:cohom} (see also \cite{mo:ergo}).
\end{exs}

\subsection{Comparison with amenability}

\begin{defn}[\cite{ra:virt}]\label{prop-moye} We say that the measured groupoid
$(G,C)$ is {\it proper} if there is a Borel system $\nu$ of {\it
probability measures} for $r : G \ra G^{(0)}$ such that $\gamma \nu^{s(\ga)} = 
\nu^{r(\ga)}$ almost everywhere. We say that $(G,C)$ is {\it amenable}
if there exists a net $(\nu_n)$ of Borel systems of probability measures such that
the net of functions $\ga \mapsto \| \ga \nu_{n}^{s(\ga)} -\nu_{n}^{r(\ga)}\|_1$
goes to $0$ in $L^\infty(G,C)$ endowed with the weak*-topology. 
\end{defn}

The notion of proper measured groupoid is the analogue in ergodic theory of the notion of compact group.
Amenable actions where introduced by Zimmer in \cite{zi:amen}. For a general
study of this notion, see \cite{dr:amenable}.

\begin{rem}\label{ess-ergo} Recall that an ergodic measured groupoid $(G,C)$
is called {\it essentially transitive} if some equivalence class $[x]$ of
units is $r(C)$-conull. Otherwise, one says that $(G,C)$ is {\it properly ergodic}.

A proper measured groupoid $(G,C)$ has a good orbit space, in the sense
that the space $G^{(0)}/G$ of equivalence classes $[x]$ is countably separated
(see \cite[Lemma 2.1.3]{dr:amenable} for instance). If in addition $(G,C)$ is
ergodic, the measure class $r(C)$ is supported by an orbit, that is, $(G,C)$ is
essentially transitive. 

In fact, the study of the structure of  amenable measured groupoids carried out
in \cite[Chapter 5]{dr:amenable} can be adapted straightforwardly to show that an
ergodic measured groupoid is proper if and only if it is essentially
transitive and for a. e. $x \in G^{(0)}$ the isotropy group $G(x)$ is compact.
\end{rem}

For the proof of the next proposition, we need to introduce an equivalent
definition of amenability for the measured groupoid $(G,C)$.
We assume, without loss of generality, that
a Haar  system $\lambda$ and a quasi-invariant measure $\mu$ are given such
$\mu\circ\lambda\in C$. For $x\in G^{(0)}$, let $\h(x) = L^2(G^x, \lambda^x)$. There
is a unique Borel structure on this bundle of Hilbert spaces such that every Borel
function $f$ on $G$, with $\int |f|^2d\lambda^x < +\infty$ for all $x$, is a Borel
section. The left translation $\hbox{Reg}(\ga)$ defined by
$\D\big (\hbox{Reg}(\ga)f\big)(\gamma_1) = f(\ga^{-1}\ga_1)$
for $\ga_1\in G^{r(\ga)}$ gives an isometry from $L^2(G^{s(\gamma)},
\lambda^{s(\gamma)})$ onto $L^2(G^{r(\gamma)},\lambda^{r(\gamma)})$.
We get a representation, independent from the choice of $(\lambda,\mu)$ with
$\mu\circ\lambda \in C$, 
up to equivalence, that is called the {\it left regular representation of}
$(G,C)$. It is proved in \cite[Theorem 6.1.4]{dr:amenable} that the
measured groupoid $(G,C)$ is  amenable if and only if its regular representation
weakly contains the trivial representation.

\begin{prop}\label{proper-am} Let $(G,C)$ be an ergodic measured groupoid. The
following conditions are equivalent:
\begin{itemize}
\item[(i)] $(G,C)$ is proper;
\item[(ii)] $(G,C)$ is amenable and has property $T$.
\end{itemize}
\end{prop}

\begin{proof} (i) $\Rightarrow$ (ii) Assume that $(G,C)$ is proper.
Obviously, $(G,C)$ is amenable. Let us
show that $(G,C)$ has property $T$. Passing to an inessential reduction we may assume
that $\gamma \nu^{s(\ga)} = \nu^{r(\ga)}$ everywhere in Definition \ref{prop-moye}, and also
that there exists a Haar system $\{\lambda^x : x \in G^{(0)}\}$
and a quasi-invariant probability measure $\mu$ on $G^{(0)}$ with
$\mu\circ \lambda \in C$.
Let us choose a non-negative Borel function $f$ on $G$ such that $\lambda^x(f) =1$
for every $x\in G^{(0)}$ and set
$$g(\ga) = \int f(\ga_{1}^{-1}\gamma) d\nu^{r(\gamma)}(\gamma_1).$$
This function is invariant by left translations and the Fubini theorem gives
$\lambda^x(g) = 1$ for all $x\in G^{(0)}$.
Let $(G^{(0)}*\h, L)$ be a representation almost having unit invariant sections.
In particular there exists a unit section $\xi$
such that
\begin{equation}\label{prov}
\int \|\xi\circ r(\ga) - L(\ga)\xi\circ s(\ga)\| g(\ga) d\mu\circ \lambda(\ga)
<1.
\end{equation}
We set $\D \eta(x) = \int L(\ga) \xi\circ s(\ga) g(\ga) d\lambda^x(\gamma)
\in \h(x)$. First,  $\eta$ is an invariant section, since
\begin{align*}
L(\ga)\eta\circ s(\ga) & = \int L(\ga \ga_1) \xi\circ s(\ga_1)
g(\ga_1) d\lambda^{s(\ga)}(\gamma_1)\\
&= \int L(\ga_1) \xi\circ s(\ga_1)
g(\ga^{-1}\ga_1) d\lambda^{r(\ga)}(\gamma_1) = \eta\circ r(\ga).
\end{align*}
Let us check that $\eta$ is non-zero. We have
$$1-\l \xi(x), \eta(x) \r  = \int \Big\l \xi(x), \xi(x) - 
L(\ga) \xi\circ s(\ga)\Big\r g(\ga) d\lambda^x(\gamma)$$
and therefore
\begin{align*}
\Big | 1 - \int \l \xi(x), &\eta(x) \r d\mu(x) \Big | \\ &=
\Big |\int_G \Big\l \xi\circ r(\ga), \Big (\xi\circ r(\ga) - 
L(\ga) \xi\circ s(\ga)\Big ) g(\ga)\Big\r d\mu\circ\lambda(\gamma)\Big |Ê\\
&\leq \int_G \|\xi\circ r(\ga) - L(\ga)\xi\circ s(\ga)\| g(\ga) d\mu\circ
\lambda(\ga) <1,
\end{align*}
by the Cauchy-Schwarz inequality and (\ref{prov}).
Thus, we get $\eta \not= 0$. Since $\eta$ is
invariant, using the ergodicity assumption, we
see that $\eta$ is a.e. equal to a non-zero constant function.

(ii) $\Rightarrow$ (i)  If $(G,C)$ is 
amenable  and has property $T$, then the regular representation
contains an invariant section, that is, there exists a Borel function
$h$ on $G$ such that $\int |h|^2 d\lambda^x = 1$ for all $x$, and 
$h(\gamma\ga_1) = h(\ga_1)$ a.e. on $G^{(2)}$. The measures $\nu^x$
 of density $|h|^2$ with respect to  $\lambda^x$ give an invariant Borel system
of probability measures for $r$, thus showing that $(G,C)$ is proper.
\end{proof}

\subsection{Cohomological characterization of property T}
The proof given by Guichardet for locally compact groups (\cite{gui:coho},
\cite[page 48]{hv}) can easily
be adapted  to our context, to show the following result. 

\begin{thm}\label{cara-suf} Let $(G,C)$ be an ergodic measured groupoid such that
for every representation $(G^{(0)}*\h, L)$ we have
$H^1((G,C), L) = 0$. Then $(G,C)$ has property $T$.
\end{thm}

\begin{proof} Let $(G^{(0)}*\h, L)$ be a representation of $(G,C)$ almost
having unit invariant sections, and assume that it has not non-zero invariant
sections. Let $\beta$ be the map from
$S((G^{(0)},r(C)),\h)$ to $B^1\big((G,C),L\big)$ defined by
$$\beta(\xi)(\ga) = \xi\circ r(\gamma) - L(\ga)\xi\circ s(\gamma).$$
This map is linear, its
range is $B^1\big((G,C),L\big)$, and it is injective by the above hypothesis.
Let us show that $\beta$ is continuous: if $(\xi_n)$ is a sequence in
$S((G^{(0)},r(C)),\h)$ going to $0$ almost everywhere with respect to $r(C)$,
obviously
$\D\big ((\beta(\xi_n)\big )$ also goes to $0$ almost everywhere with respect to $C$,
and the conclusion follows from the characterization (iii) of the topology
of the space of sections in Proposition \ref{topolo}. 

Now, since $\beta$ is a linear bijective continuous map between two
 metrizable complete topological vector spaces (Proposition \ref {metr-comp}),
its inverse
$\beta^{-1}$ is automatically continuous (see \cite[page 186]{ko} for instance).
But this is impossible, because there exists a sequence $(\xi_n)$
in $S_1((G^{(0)},r(C)),\h)$ such that $\lim_{n\ra +\infty} \beta(\xi_n) = 0$.
\end{proof}

\begin{rem} The proof actually shows that it is enough to assume that,
for every representation $(G^{(0)}*\h, L)$, the space 
$H^1((G,C), L)$ is Hausdorff. 
\end{rem}

\begin{rem} If $(G,C)$ is an amenable groupoid which is not proper, 
  then $H^1((G,C), \hbox{Reg})$ is not Hausdorff. This follows from the fact 
that the regular representation weakly contains the trivial representation without
containing it. As in the proof of Theorem \ref{cara-suf}, the corresponding map
$\beta^{-1}$ cannot be continuous, and therefore $B^1((G,C), \hbox{Reg})$
is not closed.
\end{rem}

The converse of Theorem \ref{cara-suf} is more difficult to
establish. In the case of locally compact groups, it is a result of Delorme
\cite[Th\'eor\`eme V.1]{del}. We follow the approach of \cite[Section 4.b]{hv}
suggested by Serre. The following lemma is an easy observation (see \cite[p.
49]{hv}).

\begin{lem}\label{real-comp} Let $(G,C)$ be an ergodic measured groupoid. The
following conditions are equivalent:
\begin{itemize}
\item[(i)] $H^1((G,C),L) = 0$ for every representation $L$ in a complex Hilbert
bundle;
\item[(ii)] $H^1((G,C),L) = 0$ for every representation $L$ in a real Hilbert
bundle.
\end{itemize}
\end{lem}

\begin{thm}\label{cara-nec} Let $(G,C)$ be an ergodic $r$-discrete
 measured groupoid with property $T$. Then for every
representation $(G^{(0)}*\h, L)$ of $(G,C)$ we have $H^1((G,C), L) = 0$.
In other words, every  action of $(G,C)$ by affine isometries on an Hilbert space
has an invariant section.
\end{thm}

\begin{proof} Let $(G^{(0)}*\h, L)$ be a representation of $(G,C)$.
For simplicity, and since the groupoid is ergodic, we may assume
that the Hilbert bundle is trivial, with fibre $\ck$. Also,
by the previous lemma, we can take $\ck$ to be a real Hilbert space. We shall need
the following facts, proved in \cite{hv}.

\begin{prop}[Propositions 13 and 14, p. 50-51, \cite{hv}]\label{rapp} Let $\ck$ be an
affine real Hilbert space, and $t$ a real number $>0$. 
\begin{enumerate}
\item There exists a
unique pair $(\ck_t, \phi_t)$ where $\ck_t$ is a real Hilbert space, and $\phi_t :
v\mapsto v_t$ is a continuous map from $\ck$ to $\ck_t$, such 
\begin{equation}\label{prod-scal}
\l v_t,w_t\r = \exp(-t\|v-w\|^2)
\end{equation}
for every $v,w\in \ck$, and such that $\ck_t$ is the closed vector space generated
by $\phi_t(\ck)$.
\item Given any sequence $(v_n)_{n\geq 1}$ of elements in $\ck$ such that
$\D\lim_{n\ra +\infty}
\|v_n -w\| =\infty$ for all $w\in \ck$, then the sequence $((v_n)_t)_{n\geq 1}$ goes
to
$0$ weakly in the Hilbert space $\ck_t$.
\end{enumerate}
\end{prop}

Let us come back to the proof of Theorem \ref{cara-nec}. Let $b : G^{(0)}
\ra \ck$ be a $L$-cocycle et let $\alpha$ be the homomorphism from $G$
into the group of affine isometries of $\ck$ defined by 
$\alpha(\gamma) v = L(\ga)v + b(\gamma)$ for $\ga \in G$ and $v\in \ck$.
Denote here by $\ck_n$ the Hilbert space $\ck_{1/n}$ introduced  for $t = 1/n$
in the previous proposition. The unicity of the construction shows the existence
of a representation $L_n$ into $\ck_n$ defined by
$L_n(\gamma)v_n = (\alpha(\gamma)v)_n $ for $\ga \in G$ and $v\in \ck$, where
here also we write $w_n$ instead of $w_{1/n}$ for $w\in \ck$.

We observe first that the representation $\bigoplus_{n\geq 1} L_n$
almost has unit invariant sections. Indeed, take an element $\xi\in \ck$.
Then
$$\|L_n(\gamma) \xi_n - \xi_n\|^2 = 2\Big(1-\exp(-\frac{1}{n}\|\alpha(\gamma)
\xi - \xi\|^2)\Big)$$
 by formula (\ref{prod-scal}). Hence, the sequence of functions
$(\|L_n(\cdot)\xi_n - \xi_n\|)$ goes to zero pointwise. If we identify $\xi_n$ with
the constant section taking value $\xi_n$, the conclusion follows from
characterization (ii) in  Lemma \ref{equi-sect}.

Now, since $(G,C)$ has property $T$, there exists a unit invariant section for
the representation $\bigoplus_{n\geq 1} L_n$, and at least one of its components
is non-zero. Hence, there exists an integer $m$ such that the representation
$L_m$ has a non-zero invariant section $x\mapsto\eta(x) \in \ck_m$ . 

We choose a  countable subset $\{e_k : k\in \N\}$ of the vector space
generated by $\phi_m(\ck)$ that is dense in the unit ball of $\ck_m$. Let us set
$$E_{k,n} = \{x\in G^{(0)} : |\l\eta(x), e_k\r | \geq \frac{1}{n} \}.$$
Since $x \mapsto\|\eta(x)\|$ is a non-zero constant function, there exists a pair
$(k,n)$ such that $\mu(E_{k,n}) >0$, where $\mu$ is a fixed measure in $r(C)$.
The vector $e_k$ can be written $e_k = \sum_{i=1}^p a_i (v_{i})_m$
with $a_i \in \R$ and $v_i \in \ck$. For at least one of these elements
$v_i$, there is a constant $c >0$ with $\mu(|\l \eta(\cdot), (v_{i})_m\r | \geq c)
>0$. 

Therefore, we have proved the existence of an element $v\in \ck$ and of a
constant $c>0$ such that the set
$$E:= \{x\in G^{(0)} : |\l \eta(x), v_m\r | \geq c \}$$
has a positive measure. To end the proof, it suffices to check
that for every $x\in E$, we have $\sup\{\|b(\ga)\| : \ga \in G_{E}^x \} <+\infty$
and then to use Theorem \ref{triv-cocy}. Assume on the contrary that there exists
$x\in E$ and a sequence $(\ga_n)$ in $G_{E}^x$ such that $\lim_{n\ra +\infty}
\| b(\ga_n)\| = +\infty$. We have also $\lim_{n\ra +\infty} \|\alpha(\gamma_n)v\| =
+\infty$, so that $\lim_{n\ra +\infty} \l \big(\alpha(\gamma_n)v\big)_m, \eta(x) \r
= 0$ by 
 Proposition \ref{rapp} (ii). On the other hand we have
\begin{align*}
 \Big\l \big(\alpha(\gamma_n)v\big)_m, \eta(x) \Big\r &= \Big\l L_m(\ga_n)v_m, \eta(x)\Big\r\\
&= \Big\l v_m, L_m(\ga_n)^{-1} \eta\circ r(\ga_n) \Big\r \\
& = \Big\l v_m, \eta\circ s(\ga_n) \Big\r ,
\end{align*}
since $\eta$ is an invariant section for $L_m$.
Hence we get $\lim_{n\ra +\infty} \l v_m, \eta\circ s(\ga_n) \r = 0$
which contradicts the fact that, since $s(\ga_n) \in E$, we have $| \l v_m,
\eta\circ s(\ga_n) \r | \geq c >0$ for every $n$.
\end{proof}

\section{Applications}

In this section, we give several examples of applications of Theorems
\ref{cara-suf} and \ref{cara-nec}. As mentioned in the introduction, some of these results are already known
 in case of group actions, with other proofs and sometimes in weaker forms. Some
others seem to be common folklore results that we have been unable to locate in the
litterature.

\subsection{Invariance by similarities} The invariance of property $T$
by similarity is well known for ergodic discrete measured equivalence relations
(see \cite{mo:ergo}). In this context, equivalent formulations of the notion
of similarity are given in \cite[Theorem 3]{fm1}. In particular, a more familiar way
to express similarity is the following: two
ergodic discrete measured equivalence relations $(R_1,\mu_1)$ and $(R_2,\mu_2)$
are similar if and only if there are Borel sets of positive measure
$E_1,E_2$ in $R_1^{(0)}$ and $R_2^{(0)}$ respectively such that the reduced
equivalent relations $(R_1|_{E_1},\mu_1|_{E_1})$ and  $(R_2|_{E_2},\mu_2|_{E_2})$ 
are isomorphic.

As an immediate consequence of Proposition \ref{stable},
and Theorems \ref{cara-suf} and \ref{cara-nec}, we get:

\begin{prop}\label{invsim} For ergodic $r$-discrete measured groupoids, property $T$
is invariant by similarity.
\end{prop}

In fact, we get in the same way the following result.

\begin{prop} Let $(G_i, C_i)$, $i=1,2$, be two ergodic $r$-discrete measured
groupoids and let $\phi_1 : (G_1,C_1) \ra (G_2,C_2)$, $\phi_2 : (G_2,C_2) \ra
(G_1,C_1)$ be two homomorphisms such that $\phi_1\circ\phi_2$ is similar to the
identity homomorphism of $G_2$. If $(G_1,C_1)$ is a Kazhdan measured groupoid, so
is $(G_2,C_2)$.
\end{prop}

For a transitive action of a locally compact group, one can show the stability of
property $T$, without discreteness assumptions.

\begin{thm}\label{invgroup} Let $G$ be a locally compact group and $H$ a closed
subgroup. The homogeneous $G$-space $G/H$ has property $T$ if and only if the group
$H$ has property $T$.
\end{thm}

\begin{proof} Here $G/H$ is endowed with the mesure class $C$ of any quasi-invariant
measure $\mu$ on $G/H$, and since there is a unique possible choice for $C$
we omit to mention it.  Let us recall that there is a bijective correspondence
between
$\hbox{Rep}(G/H\croi\, G)$ and $\hbox{Rep}(H)$ (see \cite[Proposition
4.2.13]{zi:ergo}). To a strict representation $L$ of $G/H\croi\, G$
in a Hilbert space $\ck$, we associate the representation $\pi_L : h\mapsto
L(\dot{e},h)$ of $H$, where $\dot{e}$ is the class of the unit element of $G$. 

To describe the inverse
map, we choose a Borel section $\st : G/H \ra G$ such that $\st(\dot{e}) = e$
and such that $\st(K)$ is relatively compact for every compact subset $K$ of $G/H$.
We set $\delta(\dot x, g) = \st(\dot x)^{-1} g \st(g^{-1} \dot x)$. It is
straightforward to check that $\delta$ is a Borel homomorphism from the groupoid
$G/H\croi\, G$ into the group $H$. Then to every $\pi \in \hbox{Rep}(H)$ we associate
$L_\pi =\pi\circ\delta \in \hbox{Rep}(G/H\croi\, G)$.
This map $\pi \mapsto L_\pi$ is the inverse of $L \mapsto \pi_L$ (up to
equivalence). It sends the trivial representation
$\pi_0$ of $H$ onto the trivial representation $L_0$ of $G/H\croi\, G$, and it is
continuous in $\pi_0$. Indeed, let $f$ be a non-negative function on $G/H\croi\, G$
with compact support $K$ and $a := \int f(\dot x,g) d\mu(\dot x) d\lambda(g) >0$,
and let us associate to $f$ the following weak*-neighbourhood of $0$:
$$\mathcal{V} = \{\varphi \in L^\infty(G/H\croi\, G) : \Big|\int \varphi(\dot x,g)
f(\dot x,g) d\mu(\dot x) d\lambda(g)\Big | \leq 1\}.$$
Let $K'$ be a compact subset of $H$ containing $\delta(K)$, and denote by
$V$ the neighbourhood of $\pi_0$ in $\hbox{Rep}(H)$ formed by the representations
$\pi$ for which there exists a unit vector $\xi$ with $\D\sup_{g\in K'} \|\pi(g) \xi
- \xi\| \leq a^{-1}$. Then, if $\tilde{\xi}$ denotes the constant function on $G/H$
with value $\xi$, we clearly have
$$\int f(\dot x,g) \|L_\pi(\dot x,g)\tilde{\xi}(g^{-1}\dot x) -
\tilde{\xi}(\dot x)\|d\mu(\dot x) d\lambda(g)\Big | \leq 1\},$$
and $L_\pi$ belongs to the neighbourhood $W(L_0, \mathcal{V})$
of $L_0$ in $\hbox{Rep}(G/H \croi\, G)$.

To prove Theorem \ref{invgroup} assume first that the $G$-space $G/H$ has property
$T$. Let
$\pi$ be a representation of $H$ almost having invariant vectors, that is, such that
$\pi_0$ belongs to the closure of $\pi$. Then $L_0$ belongs to the closure of
$L_\pi$, in other words, $L_\pi$ almost has invariant unit sections. It follows that
$L_\pi$ has an invariant unit section $\xi$, and $\pi$ has a non-zero invariant
vector by \cite[Proposition 4.2.19]{zi:ergo}.

The converse is a particular case of Corollary \ref{induce}
 proved in the next subsection, since the $G$-space
$G/H$ is induced by the action of $H$ on a space reduced to a point.
\end{proof}

\begin{cor} Let $G$ be a locally compact group, $H$ a lattice in $G$ and
$\Gamma$ a discrete infinite Kazhdan subgroup of $G$ such that the $\Gamma$-space
$G/H$ is ergodic. Then the $H$-space $G/\Gamma$ has property $T$.
\end{cor}

\begin{proof} By \cite[Proposition 2.4]{zi:cohom} (see also Corollary \ref{zimmer}
below) we know that $(H\backslash G) \croi\, \Gamma$ is a Kazhdan groupoid,
and we conclude thanks to Theorem \ref{invsim}, after having observed that
$(H\backslash G) \croi\, \Gamma$ is similar to $H\,\lcroi (G/\Gamma)$.
\end{proof}

\subsection{Mackey range}

Let us first recall a few facts on this important invariant of a cocycle.
For simplicity, we limit ourself here to the case of group actions.
In this case more details can be found in \cite[Section 4.2]{zi:ergo}. For the
general case of measured groupoids we refer to \cite[Section 7]{ra:virt}.

For the reader's convenience, we shall try to keep the notations of Zimmer's book
\cite{zi:ergo} as much as possible, and in particular in this subsection we shall
consider right actions. Let $G$ and $H$ be two locally compact spaces,
$(S,\mu)$ a measured
$G$-space and $\alpha : S\croi\, G \ra H$ a homomorphism (or cocycle) in the sense of
Definition \ref{hom}. We endow $S\times H$ with the class of the measure $\nu =\mu\times
\lambda$. The {\it skew-product action} is the action of $G$ on $(S\times H, \nu)$
defined by $(s,h)g = (sg , h\alpha(s,g))$. This measured $G$-space is usually denoted
by $S\times_\alpha H$. There is also an action of $H$ on $S\times H$ defined
by $(s,h)h_1 = (s,h_{1}^{-1}h)$, which commutes with the $G$-action. The {\it Mackey
range} or {\it range closure} of $\alpha$ is the action of $H$ on the standard
quotient $X$ of $S\times_\alpha H$ with respect to the $G$-action (i.e. the mea\-sured space
that realizes the ergodic decomposition of $\nu$ with respect to the $G$-action). In
case the space $(X\times_\alpha G)/G$ of $G$-orbits is countably separated, note that
$X$ and $(S\times_\alpha H)/G$ are isomorphic as measured $H$-spaces. We shall denote
by $p$ the natural map from $S\times H$ on $X$ and by $[\dot{\nu}]$ the class of the image
under $p$ of any probability measure equivalent to $\nu$. Observe that $p$ is measure preserving
and $H$-equivariant.

Let us also
recall that the Mackey range of an ergodic action is ergodic.

\begin{thm}\label{invMackey0} Let $G$ and $H$ be two locally compact groups. Let $(S,\mu)$ be an ergodic (right) $G$-space
such that for every representation $L$ of the measured groupoid $S\croi\, G$
we have $H^1((S\croi\, G, [\mu]), L) = 0$. Then for every representation $L$ of
the Mackey range $X\croi\, H$ of a homomorphism $\alpha : (S\croi\,G,\mu) \ra H$,
we still have  $H^1((S\croi\, H, [\dot{\nu}]), L) = 0$. 
  \end{thm}

\begin{proof} As said before, we consider here right group actions. We shall follow
the same steps as in the proof of
\cite[Theorem 3.3]{zi:amen} showing that the Mackey range of every amenable action
is amenable. When possible, we shall choose the same notations as in this proof, to
facilitate the comparison. Let $L$ be a representation of $X\croi\, H$ in a Hilbert
space $\ck$,
$b$ a $L$-cocycle, and $\gamma$ the corresponding affine action. We can assume
that $L$, $b$ and therefore $\gamma$ are strict on an inessential reduction, here denoted
by $X_1 * H :=\{(x,h) : x\in X_1, xh \in X_1\}$ as in \cite{zi:amen}. The crucial
technical steps are provided by \cite[Lemmas 3.5, 3.6]{zi:amen}, based on Ramsay's
results. First, there is an inessential reduction $S_0*G$ of $S\croi\, G$  and a
homomorphism $\tilde{\alpha} : S\croi\, G \ra H$ such that for all $g\in G$, 
$\tilde{\alpha}(s,g) = \alpha(s,g)$ a.e., and such that
$(\tilde{p},\tilde{\alpha}) : S_0*G \ra X\croi\, H$
is a strict homomorphism, where $\tilde{p}(s) = p(s,e) \in X$.
Now, we replace $(\tilde{p},\tilde{\alpha})$ by a cohomologous homomorphism
taking values in $X_1*H$. By \cite[Lemma 3.6]{zi:amen}, there is a conull
Borel subset $S_1 \subset S_0$ and a Borel function $\theta : S_1 \ra H$
such that, if we set $q(s) = \tilde{p}(s)\theta(s)^{-1}$ and $\beta(s,g)
= \theta(s) \tilde{\alpha}(s,g)\theta(sg)^{-1}$, then $(q,\beta)$
is a strict homomorphism from $S_1*G$ into $X_1*H$. Observe that
$\delta :=\gamma\circ(q,\beta)$is a homomorphism from the inessential reduction
$S_1*G$ of $S\croi\, G$ into the group of affine isometries of $\ck$, corresponding
to the representation $L\circ(q,\beta)$. Therefore there exists a $\delta$-invariant
Borel map $\varphi : S\ra \ck$.

Next, we define $\psi : S_1 \times H \ra \ck$ by 
$$\psi(s,h) = \ga(q(s), h^{-1})^{-1}\varphi(s).$$
As in \cite{zi:amen}, for $(s,g) \in S_1\croi\, G$ and for all $h\in H$, we have
\begin{align*}
\psi(sg,h\beta(s,g)) &= \ga(q(sg), \beta(s,g)^{-1}h^{-1})^{-1} \varphi(sg)\\
& = \ga(q(s),h^{-1})^{-1} \ga(q(sg),\beta(s,g)^{-1})^{-1}\varphi(sg)\\
&= \ga(q(s),h^{-1})^{-1} \ga(q(s),\beta(s,g))\varphi(sg)\\
&= \ga(q(s),h^{-1})^{-1} \varphi(s) = \psi(s,h),
\end{align*}
because $q(sg) = q(s)\beta(s,g)$, since $(q,\beta)$ is a homomorphism.

Now, let $\omega(s,h) = \psi(s, h\theta(s)^{-1})$. Then for $(s,g) \in S_1*G$
and almost all $h\in H$ we have
\begin{align*}
\omega(sg, h\tilde{\alpha}(s,g))&= \psi(sg, h\tilde{\alpha}(s,g)\theta(sg)^{-1})\\
&= \psi(sg, h\theta(s)^{-1}\beta(s,g))\\
&= \psi(s,h\theta(s)^{-1}) = \omega(s,h).
\end{align*}
Therefore, $\omega$ is essentially $G$-invariant on $S\times_\alpha H$,
and hence there is a map $\sigma : X \ra \ck$ such that $\sigma(p(s,h)) =
\omega(s,h)$ a. e. To conclude the proof, one checks exactly as in \cite[page
365]{zi:amen} that for every $h\in H$,
$$\ga(s,h)\sigma(sh) = \sigma(s) \quad\hbox{a. e.,}$$
and this shows that $b$ is a $L$-coboundary.
\end{proof}

\begin{cor}\label{invMackey} Let $G$ be a discrete group, $(S,\mu)$ an ergodic
$G$-space having property $T$ and $H$ a locally compact group. The Mackey range
of every homomorphism $\alpha : (S\croi\,G,\mu) \ra H$ is a $H$-space having
property $T$.
\end{cor}

\begin{rem} This result is contained in \cite{ne} under the additional assumption
that $H$ is discrete. Nevo also proved this corollary when $H$ is discrete and
$G$ is a locally compact Kazhdan group preserving a probability measure $\mu$.
\end{rem}

\begin{cor}\label{induce} Let $H$ be a closed subgroup of a locally compact
group $G$, and let $(Y,\mu)$ be an ergodic measured $H$-space such that
$H^1((Y\croi\, H, [\mu]), L) = 0$ for every representation $L$. 
Then the induced $G$-action has property $T$.
\end{cor}

\begin{proof} For the definition of an induced action, we refer to
\cite[Definition 4.2.21]{zi:ergo} or to \cite{zi:induce}. The result follows
immediately from Theorem \ref{invMackey}, after having observed that the induced
action from $H$ to $G$ is the Mackey range of the homomorphism $(y,h)\mapsto h$
from $Y\croi\, H$ into $G$.
\end{proof}

\begin{cor}\label{skewergo} Let $G$ be a discrete group, $(S,\mu)$ an ergodic
$G$-space having property $T$, and $\alpha$ a homomorphism
from $(S\croi\, G,\mu)$ into a locally compact group $H$ such that the skew-product
$S\times_\alpha H$ is an ergodic $G$-space. Then the group $H$ has
property $T$.
\end{cor}

\begin{proof} It suffices to observe that the ergodicity of the skew-product means
that the Mackey range is the action of $H$ on a point. This action has property $T$
exactly when $H$ has property $T$.
\end{proof}

One can extend Corollary  \ref{invMackey} and prove in the same way the following
result.

\begin{thm}\label{invMackey2} Let $(G,C)$ be an ergodic $r$-discrete measured Kazhdan
groupoid and
$H$ a locally compact group. Then the Mackey range of every homomorphism from
$(G,C)$ to $H$ is a $H$-space having property $T$.
\end{thm}

\begin{rem} In \cite[Theorem 10]{zi:kazhdan}, Zimmer has proved the following 
remarkable result: let $G$ be a discrete Kazhdan group acting on a probability
space $(X,\mu)$ so as to preserve the probability measure $\mu$. Let $H$ be
a real algebraic group and $\alpha : X\croi\, G \ra H$ a cocycle. Then $\alpha$ is
cohomologous to a cocycle $\beta$ such that $\beta(X\croi\, G) \subset H_1\subset H$
where $H_1$ is an algebraic Kazhdan subgroup of $H$.  In other words, the Mackey
range of $\alpha$ is an action of $H$ induced by a Kazhdan algebraic subgroup of
$H$.  This result is to be compared to our Corollary \ref{invMackey}
\end{rem}

\subsection{Extensions and quotients} 
Let $(G,\lambda)$ be a Borel groupoid with Haar system, and let $Y$, $X$ 
be two (left) Borel $G$-spaces respectively equipped with  measures $\nu$ and $\mu$
quasi-invariant for $(G,\lambda)$. If in addition we are given 
an equivariant Borel map $p: Y \ra X$  such that $\mu$ is equivalent to a
pseudo-image of $\nu$ by $p$, we say that $((Y,\nu),(X,\mu), p)$ (or simply
$(Y,X)$) is a {\it pair of measured $G$-spaces}, or that $Y$ {\it
is an extension of} $X$. If both $G$-spaces are ergodic, we say that the pair
is ergodic.

 If $\nu$ admits the disintegration $\nu = \int \rho^x d\mu(x)$ along
$p$, where $\{\rho^x : x\in X\}$ an invariant Borel system of probability measures
for $p$ we say that $(Y,X)$ is a {\it measure preserving  pair of $G$-spaces}.

\begin{lem}\label{injpair} We keep the previous notations and we assume that $(Y,X)$
be a measure preserving ergodic pair of
$G$-spaces. Denote by $\phi$ the homomorphism $(y,\ga) \mapsto
(p(y),\ga)$ from $Y\croi\, G$ onto $X\croi\, G$. Then for every
representation
$L$ of $(X\croi\, G, [\mu])$, the map $\phi^* : H^1((X\croi\, G, [\mu]), L)
\ra  H^1((Y\croi\, G, [\nu]), L\circ\phi)$ is injective.
\end{lem}

\begin{proof}
 We can assume that $L$ is a representation  into the constant field
$\h$ with fibre $\ck$. Let $b : X\croi\, G \ra \ck$ be a $L$-cocycle, and denote
by
$\alpha$ the corresponding action of $X\croi\, G$ by affine transformations of $\ck$.
Assume that the cocycle $b\circ \phi$ is a $L\circ \phi$-coboundary.
This means that, passing if necessary to
an inessential reduction, there exists a Borel map $\xi : Y \ra \ck$ such that
\begin{equation}\label{inva}
\xi(y) = \alpha(p(y), \ga) \xi(\ga^{-1}y),
\end{equation} 
that is
\begin{equation}\label{inva2}
\xi(y) =L( p(y), \ga) \xi(\ga^{-1}y) + b(p(y), \ga)
\end{equation}  
 for $(y,\ga) \in Y\croi\, G$. We want to show that $b$ is a $L$-coboundary.

Note first that the subset
$X_1\subset X$ of all elements
$x$ such that $y\mapsto \|\xi(y)\|$ is essentially bounded on $p^{-1}(x)$
is invariant. Indeed, consider two equivalent elements $x$ and $\ga^{-1}x$
in $X$. Then we have $y\in p^{-1}(x)$ if and only if $\ga^{-1}y \in
p^{-1}(\ga^{-1}x)$, and from (\ref{inva2}) above we get
$$\|\xi(y)\| \leq \|\xi(\ga^{-1}y)\| +\|b(x, \ga)\|$$
 for every $y \in p^{-1}(x)$, hence our assertion. 
Now,  since the action is ergodic, $X_1$ is either null or conull. The main point is to prove that $X_1$
is conull.

Let us assume, on the contrary, that $X_1$ is null. We choose a ball $B$ in $\ck$,
centered to
$0$, whose radius $R$ is large enough,  so that $E = \xi^{-1}(B)$ has a positive
$\nu$-measure. Since $\int \rho^x(E) d\mu(x) >0$, there exists $\beta >0$
such that $A :=\{x\in X : \rho^x(E)\geq \beta \}$ has a positive $\mu$-measure. 
We set $E' = E \cap p^{-1}(A)$.

Let $y \in [E']$ and choose $(y,\ga) \in Y\croi\, G$ such that
 $\ga^{-1}y\in E'$. If we set $x = p(y)$, we have $\ga^{-1}x = p(\ga^{-1}y) \in A$
and 
$$\xi(y) = \alpha(x, \ga) \xi(\ga^{-1}y) \in \alpha(x, \ga)(B)$$
by (\ref{inva}).

Since $[E']$ is conull and since $p([E']) = [A]$,
there is a conull subset $A_1$ of $[A]$ such that $\rho^x([E']) = 1$ for $x\in A_1$.
In particular, we get the following fact: for $x\in A_1$ and for $\rho^x$-almost
every $y\in p^{-1}(x)$ we have 
\begin{equation}\label{inclus}
\xi(y) \in \bigcup_{\{\ga: \ga^{-1}x \in A\}}\alpha(x,\ga)(B).
\end{equation}

This implies that for almost every $x\in A_1$, there is an infinite subset
$\{(x,\ga_i) : i\in \N\}$ of $X\croi\,G$ with $\ga_{i}^{-1}x \in A$ for every $i$
and 
$\alpha(x,\ga_i)(B)\cap\alpha(x,\ga_j)(B)=\emptyset$ if $i\not = j$. Indeed, let
$x\in A_1$ for which this is not true, so that we have a maximal finite set
$\{(x,\ga_1),\dots, (x,\ga_N)\}$  of such elements. Then, for
$(x,\ga) \in X\croi\,G$ with $\ga^{-1}x\in A$, there exists $i\in \{1,\dots N\}$,
and $v, v_i \in B$ such that $\alpha(x,\ga)v = \alpha(x,\ga_i)v_i$,
that is 
$$ L(x,\ga)v +b(x,\ga) = L(x,\ga_i)v_i +b(x,\ga_i).$$
Then, for $w\in B$, we have
\begin{align*}
\alpha(x,\ga)w &= L(x,\ga)(w-v) + L(x,\ga)v + b(x,\ga)\\
&= L(x,\ga)(w-v) + L(x,\ga_i)v_i + b(x,\ga_i)
\end{align*}
and therefore
 $$\|\alpha(x,\ga)w\| \leq 3R +\|b(x,\ga_i)\|\leq 3R +\max\{\|b(x,\ga_i)\| : 1\leq i
\leq N\}.$$
 As a result, $\bigcup_{\{\ga: \ga^{-1} x \in A\}}\alpha(x,\ga)(B)$
is bounded, and therefore, by (\ref{inclus}) above, $\xi$ is essentially bounded
on $p^{-1}(y)$, so that $x \in X_1$. Since $X_1$ is supposed to be null,
our assertion follows, that is, for almost every $x \in A_1$,
there is an infinite set $\{(x,\ga_i) : i\in \N\}$ of $X\croi\,G$ with
$\ga_{i}^{-1}x \in A$ and
 $\alpha(x,\ga_i)(B) \cap \alpha(x,\ga_j)(B) = \emptyset$ if $i\not= j$.

On the other hand, let $x\in A_1$ for which there exists such an  infinite set $\{(x,\ga_i) : i\in \N\}$.
We shall show that such a $x$ cannot exist. 
We set $$D_i = \{ y \in p^{-1}(x) : \xi(y) \in \alpha(x,\ga_i)(B)\}.$$
These subsets do not intersect, and therefore we have
\begin{equation}\label{disjoint}
1\geq \rho^x(\bigcup_i D_i) = \sum\rho^x(D_i).
\end{equation}
Now, observe that
\begin{align*}
D_i &= \{y \in p^{-1}(x) : \alpha(\ga_{i}^{-1}x, \ga_{i}^{-1})\xi(y) \in B\}\\
& = \{y \in p^{-1}(x) : \xi(\ga_{i}^{-1}y) \in B\},
\end{align*}
by (\ref{inva}). It follows that 
\begin{align*}
\rho^x(D_i) &= \rho^{x}(\{\ga_i y : \xi(y)\in B\})\\
&=\ga^{-1}_i \rho^x(\xi^{-1}(B)) = \rho^{\ga^{-1}_i x}(\xi^{-1}(B)),
\end{align*}
by the invariance of $\rho$.
Since $\ga^{-1}_i x \in A$, we have 
$$\rho^{\ga^{-1}_i x}(\xi^{-1}(B)) \geq \beta> 0.$$
 Because $\rho^x$ is a probability measure,
the inequality (\ref{disjoint}) leads to a contradiction.

Hence $A_1$ is null, but this is impossible since $A_1$ is also conull. As a
consequence, we get that
$X_1$ is conull. For
$x\in X_1$, we can define
$$\eta(x) = \int \xi(y) d\rho^x(y),$$
since $\xi$ is essentially bounded on $p^{-1}(x)$.
We have
\begin{align*}
\alpha(x,\ga) \eta(\ga^{-1}x) & = \alpha(x,\ga)\int \xi(y) d\rho^{\ga^{-1}x}(y)\\
&= \alpha(x,\ga)\int\xi(\ga^{-1}y) d\rho^x(y)\\
&= \eta(x),
\end{align*}
by (\ref{inva}).
This shows that $\eta$ is invariant under the affine action, and therefore
$b$ is a coboundary.
\end{proof}

This lemma suggests the following definition.

\begin{defn}\label{kazhdanpair} Let $G$ be an $r$-discrete groupoid. We say
that an ergodic pair $(Y,X)$ of measured $G$-spaces is a {\it Kazhdan pair}
(or has {\it property} $T$) if for every representation $L$ of $X\croi\, G$, the
map
$$\phi^* : H^1((X\croi\, G, [\mu]), L) \ra 
H^1((Y\croi\, G, [\nu]), L\circ\phi)$$
 is injective.
\end{defn}

\begin{rem} By lemma \ref{injpair} a measure preserving ergodic pair of $G$-spaces
has property $T$. Furthermore an $r$-discrete ergodic measured  groupoid
$(G,\mu)$ is Kazhdan if an only if the pair
$((G,\mu\circ\lambda),(G^{(0)},\mu), r)$ is Kazhdan. This follows from the
fact that the groupoid
$G\croi\,G$ has trivial cohomology groups with respect to its
representations since it is proper. Therefore the injectivity of $\phi^* :
H^1((G, [\mu]), L) \ra  H^1((G\croi\, G, [\mu\circ\lambda]),
L\circ\phi)$ for every representation
$L$ of $G$ is equivalent to property $T$.

On the other hand, for a discrete Kazhdan group $G$, an ergodic pair $(Y,X)$ of measured $G$-spaces
is not necessarily a Kazhdan pair. For instance, if $G$ has a subgroup $H$ without property $T$, 
the pair $(G,G/H)$ cannot have property $T$. Otherwise, since the $G$-space $G$ has trivial cohomology, 
the $G$-space $G/H$ would have property $H$, that is, $H$ would be a Kazhdan group, a contradiction.
\end{rem}

\begin{thm}\label{quotient} Let $(G,\lambda)$ be a Borel groupoid with Haar
system and consider a measure preserving ergodic pair $((Y,\nu),(X,\mu),p)$ of
$G$-spaces.
\begin{itemize}
\item[(i)] If the $G$-space $X$ has property $T$ then the $G$-space $Y$
also has property $T$.
\item[(ii)] Assume that $G$ is $r$-discrete and that the $G$-space
$Y$ has property $T$. Then the $G$-space $X$ also has  property $T$.
\end{itemize}
\end{thm}

\begin{proof} Let us first recall that given a representation $L$
of $(Y\croi\, G,[\nu])$ in a Hilbert bundle $\h$, we can associate an {\it induced
representation} $U_L$ of $(X\croi\, G,[\mu])$. Without loss of generality we assume
that the bundle $\h$ is trivial, with fibre $\ck$. We denote by $p_*(\h)$ the Hilbert
bundle on $X$ with fibres $p_*(\h)(x) = L^2(\rho^x)\otimes \ck$, and $U_L$
is the representation of $X\croi\, G$ into $p_*(\h)$ defined by
$$\Big(U_L(x,\ga)\xi\Big)(y) = L(y,\ga)\xi(\ga^{-1}y)$$
for $\xi \in L^2(\rho^{\ga^{-1}x})\otimes \ck$ and $y\in p^{-1}(x)$,
so that $U_L(x,\ga)$ is an isometry from $L^2(\rho^{\ga^{-1}x})\otimes \ck$
onto $L^2(\rho^x)\otimes \ck$. For more details we refer to \cite{ra:nontran}.

Let us show that if $L$ almost contains invariant unit sections, this is also
the case for $U_L$. First, let $\xi: Y\ra \ck$ be a unit Borel map, and set
$\tilde{\xi}_x = \xi|_{p^{-1}(x)}$. Since 
$$\| \tilde{\xi}_x \|^2 = \int \|\xi(y)\|^2 d\rho^x(x) =1,$$
we see that $\tilde{\xi} : x\mapsto \tilde{\xi}_x $ is a Borel unit section of the
bundle $p_*(\h)$. Moreover, we have
$$\Big(U_L(x,\ga)\tilde{\xi}_{\ga^{-1}x}\Big)(y) = L(y,\ga)\xi(\ga^{-1}y)$$
for $y\in p^{-1}(x)$ and
$$\|U_L(x,\ga)\tilde{\xi}_{\ga^{-1}x} - \tilde{\xi}_x \|^2
= \int \| L(y,\ga) \xi(\ga^{-1}y) - \xi(y) \|^2 d\rho^x(y).$$
Second, let $(\xi_n)$ be a sequence of unit maps $\xi_n : Y \ra \ck$ such that
$\lim_n \|c_L(\xi_n)\|$ $= 0$ in $L^\infty(Y\croi\, G)$ endowed with the
weak*-topology. Then $\lim_n \| c_{U_L}(\tilde{\xi}_n)\| = 0$ in
$L^\infty(X\croi\, G)$ endowed with the same topology. Indeed, given
$f: X\croi\, G \ra \R^+$ such that $\int f d\mu\circ \lambda =1$, this is a
consequence of the following inequalities (where $\nu = \mu\circ\rho$):
\begin{align*}
\Big(\int f(x,\ga) &\|U_L(x,\ga)\tilde{\xi}_{\ga^{-1}x} - \tilde{\xi}_x \|
d\mu\circ\lambda(x,\ga)\Big)^2 \\
&\leq
\int f(x,\ga) \|U_L(x,\ga)\tilde{\xi}_{\ga^{-1}x} - \tilde{\xi}_x \|^2
d\mu\circ\lambda(x,\ga)\\
&= \int f(x,\ga) \Big[\int \| L(y,\ga) \xi(\ga^{-1}y) - \xi(y) \|^2 d\rho^x(y)\Big]
d\mu\circ\lambda(x,\ga)\\
&\leq 2 \int f(p(y),\ga) \| L(y,\ga) \xi(\ga^{-1}y) - \xi(y) \| d\nu\circ
\lambda(y).
\end{align*}

To prove (i), assume that $L$ almost contains invariant unit sections
and that $X\croi\, G$ has property $T$. It follows
from what has been said above, that $U_L$ has an invariant unit section $\eta : x \mapsto \eta_x
\in L^2(\rho^x)\otimes \ck$. Let us define $\xi : y \ra \ck$ by $\xi(y)
= \eta_{p(y)}(y)$. We have $U_L(x,\ga)\eta_{\ga^{-1}x} = \eta_x$
for $\mu\circ\lambda$-almost every $(x,\ga)\in X\croi\, G$, and therefore
$$L(y,\ga)\xi(\ga^{-1}y) = \xi(y), \quad \rho^x-a.e.$$
for $\mu\circ\lambda$-almost every $(x,\ga)$. By the Fubini theorem, we see that 
$L(y,\ga)\xi(\ga^{-1}y)$ $= \xi(y)$ for $\nu\circ\lambda$ almost every $(y,\ga) \in
Y\croi\, G$, and the proof is achieved, since $\xi$ is non-zero.

(ii) is a consequence of Lemma \ref{injpair}.
\end{proof}

\begin{cor}[Proposition 2.4, \cite{zi:cohom}]\label{zimmer} Let $G$ be a locally 
compact group and $Y$
an ergodic $G$-space with finite invariant mesure.
\begin{itemize} 
\item[(i)] If $G$ has property $T$, so does the $G$-action on $Y$;
\item[(ii)] If $G$ is a countable discrete group, and if the $G$-action has property
$T$, then $G$ has property $T$.
\end{itemize}
\end{cor}

\begin{proof} In the previous theorem, we take $X$ reduced to a single element.
\end{proof}

Note that in \cite[Proposition 2.4]{zi:cohom}, the additional assumption
that the action is mixing was needed to get (ii).

\begin{rem}\label{ME} One way to understand Gromov's measure equivalence of
countable groups is the following (see \cite[Theorem 3.3]{fu2} and \cite[Lemma 2.2]
{fu1}): two countable groups $\Gamma_1$ and $\Gamma_2$ are $ME$ (measure equivalent)
if and only if there exists ergodic essentially free finite measure preserving
actions $(X_1, \Gamma_1,\mu_1)$ and $(X_2, \Gamma_2,\mu_2)$ such that the
measured groupoids $(X_i\croi\,\Gamma_i, \mu_i)$, $i=1,2$, are similar.
Then it follows from Proposition \ref{invsim} and Corollary \ref{zimmer} 
that property $T$ is a $ME$ invariant \cite[Corollary 1.4]{fu1}.
\end{rem}

\begin{thm} Let $(G,\mu)$ be a $r$-discrete ergodic Kazhdan groupoid. Then the
associated discrete measured equivalence relation $(R,\mu)$ also has property $T$.
\end{thm}

\begin{proof} Denote by $\phi$ the homomorphism $\ga \mapsto (r(\ga),s(\ga))$
from $G$ onto $R$.  Let $L$ be a representation of $R$ into the constant field
$\h$ with fibre $\ck$ and let $b : R \ra \ck$ be a $L$-cocycle .
Since $(G,\mu)$ has property $T$, the cocycle $b\circ \phi$ is a $L\circ \phi$-coboundary.
This means that, passing if necessary to
an inessential reduction, there exists a Borel map $\xi : G^{(0)} \ra \ck$ such that
$$\xi(r(\ga)) =L\circ\phi(\ga) \xi(s(\ga)) + b\circ \phi(\ga)$$  
 for $\ga\in G$. It follows that $b$ is a coboundary for the representation $L$.
\end{proof}

Using the same techniques, it is also possible to deal with semi-direct products.

\subsection{Property $T$ and negative definite functions} As for locally
compact groups, property
$T$ for an ergodic $r$-discrete groupoid can be detected by the behaviour of its 
conditionally negative definite functions. Let us first recall some definitions and
facts.

\begin{defn}\label{cond-neg}  A {\it real conditionally negative
definite  function on a Borel grou\-poid} $G$ is a Borel function $\psi : G \ra \R$
such that we have 
\begin{enumerate}
\item $\psi(x) = 0$ for every $x\in G^{(0)}$;
\item $\psi(\gamma) = \psi(\ga^{-1})$ for  every $\ga \in G$;
\item for every $x\in G^{(0)}$, every integer $n\geq 2$, every $\ga_1,
\dots,\ga_n \in G^x$ and every real numbers $\lambda_1,\dots,\lambda_n$ with
$\sum_{i=1}^n \lambda_i = 0$, then
$$\sum_{i=1}^n \sum_{j=1}^n \lambda_i\lambda_j \psi(\gamma_{i}^{-1}\ga_{j}) \leq 0.$$
\end{enumerate}

A {\it real conditionally negative
definite  function on a measured grou\-poid} $(G,C)$ is a Borel function $\psi : G \ra \R$
such that the restriction of $\psi$ to some inessential reduction satisfies
the above three conditions.
\end{defn}

Note that $\psi$ must be a  non-negative function.

\begin{rem} Let $b$ be a $L$-cocycle for a representation $L$ of $(G,C)$. Exactly as for groups,
the map $\ga \ra \|b(\ga)\|^2$ is a conditionally negative definite  function. Given
$\ga_1,\ga_2\in G^x$, this follows from the equalities
\begin{align*}
b(\ga_{1}^{-1} \ga_2) & = b(\ga_1^{-1}) + L(\ga_{1}^{-1})b(\ga_2)\\
& = -L(\ga_{1}^{-1})b(\ga_1) + L(\ga_{1}^{-1})b(\ga_2),
\end{align*}
so that
$$ \|b(\ga_{1}^{-1} \ga_2)\|^2 = \|b(\ga_1)\|^2 + \|b(\ga_2)\|^2 -\l b(\ga_1),
b(\ga_2)\r 
 -\l b(\ga_2), b(\ga_1)\r .$$
\end{rem}

Conversely we have
\begin{prop}\label{cond-neg1} Let $G$ be a $r$-discrete Borel groupoid and let $\psi$
be a conditionally negative definite  function. There exists a representation in a
real Hilbert bundle, and a $L$-cocycle $b$, such that $\psi(\ga) = \|b(\ga)\|^2$ for
every $\ga\in G$.
\end{prop}

\begin{proof} This is essentially the proof given for groups (see \cite[Page 63]{hv}).
For $x\in G^{(0)}$, let $F(x)$ be the vector space of all maps $f : G^x \ra \R$, with finite support,
and such that $\sum_{\ga\in G^x} f(\ga) = 0$. We endow $F(x)$ with the non-negative
bilinear symmetric form
$$\l f_1, f_2\r = -\frac{1}{2} \sum_{\ga_1, \ga_2 \in G^x} f_1(\ga_1)f_2(\ga_2) \psi(\ga_{1}^{-1}\ga_2).$$
We denote by $\h(x)$ the real Hilbert space obtained by separation and completion, and by $\Lambda_x$
the map $\ga \mapsto \delta_\ga - \delta_x$ from $G^x$ to $\h(x)$ (where $\delta_\ga$ is the Dirac function at
$\gamma$ and where, by a slight abuse, we use the same notation for an element of $F(x)$
and for its image in $\h(x)$). Observe
that
$\h(x)$ is generated by the image of
$\Lambda_x$, and that
$$\|\Lambda_x(\ga_1) -\Lambda_x(\ga_2)\|^2 = \psi(\ga_1^{-1}\ga_2)$$
for $\ga_1, \ga_2 \in G^x$. Moreover, $\h(x)$ is characterized by these properties (see \cite{hv}).
It follows that for $\ga \in G$ there exists a unique orthogonal operator $L(\ga) : \h(s(\ga)) \ra \h(r(\ga))$
such that $L(\ga)(\Lambda_{s(\ga)}(\ga')) = \Lambda_{r(\ga)}(\ga\ga')$ for every
$\ga'\in G^{s(\ga)}$. Given a countable family $(\st_n)$ of Borel sections $\st_n :
x\in G^{(0}) \mapsto \st_n(x) \in G^x$, such that $\D G = \cup \st_n(G^{(0)})$,
we endow the Hilbert bundle $\h$ with the unique Borel structure making the sections
$x \mapsto \Lambda_x(\st_n(x))$ Borel. Then $L$ is a Borel representation in the
bundle $\h$, and if $b$ denotes the map $\ga \mapsto
\Lambda_{r(\ga)}(\ga)$, we immediately see that $b$ is a Borel
$L$-cocycle and that $\|b(\ga)\|^2 = \psi(\ga)$ for all $\ga \in G$.
\end{proof}

\begin{thm}\label{cond-neg-T} Let $(G,C)$ be an ergodic $r$-discrete measured
groupoid. The following conditions are equivalent:
\begin{itemize}
\item[(i)] $(G,C)$ has property $T$.
\item[(ii)] For every real conditionally negative
definite  function $\psi$, there exists a Borel subset $E\subset G^{(0)}$ of positive
measure, such that the restriction of $\psi$ to $G|_E$ is bounded.
\item[(iii)] For every real conditionally negative
definite  function $\psi$, there exists a Borel subset $E\subset G^{(0)}$ of positive
measure, such that $\D\sup_{\ga\in G_{E}^x} |\psi(\ga)|<+\infty$ for every $x\in E$.
\end{itemize}
\end{thm}

\begin{proof} (i) $\Rightarrow$ (ii) By the previous proposition, $\psi$ is of the 
form $\ga \mapsto b(\ga)$ where $b$ is a cocycle with respect to a representation
$L$ into a real Hilbert bundle $\h$. If $G$ has property $T$, there is a section
$\xi$ of $\h$ such that $b(\ga) = \xi\circ r(\ga) -  L(\ga)\xi\circ s(\ga)$
for almost every $\ga \in G$. There is a constant $c >0$ such that
$E:= \{x\in G^{(0)} : \| \xi(x)\| \leq c \}$ has positive measure. Then for
$\ga \in G|_E$, we have $\| b(\ga) \| \leq 2c$.

(ii) $\Rightarrow$ (iii) is obvious, and (iii) $\Rightarrow$ (i) follows immediately
from Proposition \ref{cond-neg1}, and Theorems \ref{triv-cocy}, \ref{cara-suf}.
\end{proof}

\begin{cor}\label{discond-neg2} Let $(G,C)$ be an ergodic $r$-discrete measured
groupoid with property $T$ and let $Z$ be a $G$-bundle of complete metric spaces over
$G^{(0)}$. We assume that each metric space $(Z(x), d_x)$ satisfies the median
inequality and that $d_x : Z(x) \times Z(x) \ra \R$ is a conditionally negative
definite  kernel. If
$Z$ has a Borel section, then there is a conull subset $U$ of
$G^{(0)}$ and a Borel section $\xi : x \in U \mapsto \xi(x) \in Z(x)$ such that
$\ga \xi\circ s(\ga) = \xi\circ r(\ga)$ for every $\ga \in G|_U$.
\end{cor}

\begin{proof} For the definition of a conditionally negative definite 
kernel, see \cite[page 62]{hv} for instance. Let $\st$ be a Borel section for $p : Z
\ra G^{(0)}$, and for $\ga \in G$ let us set $\D\psi(\ga) =
d_{r(\ga)}\big(\st\circ r(\ga),
\ga\st\circ s(\ga)\big)$. Since $\psi$ is real conditionally negative definite,
there exists a subset $E$ of $G^{(0)}$ of positive measure such that
$\sup_{\ga\in G_{E}^x} |\psi(\ga)|<+\infty$ for every $x\in E$.
The conditions of Lemma \ref{inv-sect} are fulfilled for $(G|_E, C|_E)$,
 $Z_E = p^{-1}(E)$ and $B(x) = \{\ga \st(s(\ga)): \ga\in G_{E}^x\}$, $x\in E$. 
Let $\xi : x\in E \mapsto \xi(x) \in p^{-1}(x)$ such that $\ga \xi\circ s(\ga)
= \xi\circ r(\ga)$ for $\ga\in G|_E$ (up to an inessential reduction).

 Now, using Theorem \ref{compact} and Lemma 
\ref{federer}, we may assume the existence of a Borel section $\theta : [E]
\ra r^{-1}(E)$ for $s$. Moreover, we may take $\theta(x) = x$ if $x\in E$.
For $x\in [E]$ we put
$$\tilde{\xi}(x) = \theta(x)^{-1} \xi\circ r(\theta(x)).$$
Note that $\tilde{\xi}$ extends $\xi$ and is $G|_{[E]}$-invariant. Indeed,
let $\ga \in G|_{[E]}$, and set $x=r(\ga)$, $y= s(\ga)$. Since
$\theta(x)\ga\theta(y)^{-1} \in G|_E$, we have
$$\theta(x)\ga \theta(y)^{-1} \xi\circ r(\theta(y))
= \xi\circ r(\theta(x)),$$
and therefore
$$\ga \theta(y)^{-1} \xi\circ r(\theta(y)) = \theta(x)^{-1}\xi\circ r(\theta(x)).$$
Hence, the corollary is proved with $U = [E]$, that is conull since
$(G,C)$ is ergodic.
\end{proof}

\subsection{Property $T$ and trees}

A {\it $G$-bundle of trees} is a Borel $G$-bundle $Z$ of metric spaces, such that 
each fibre $Z(x)$ is a countable tree which defines the metric $d_x$.
Then
$$K := \{(z,z') \in Z*Z  : d_{p(z)}(z,z') =1\}$$
is a symmetric Borel subset of $Z*Z$, the subset of {\it edges}, which is
$G$-invariant, that is, $(\ga z, \ga z') \in K$  for every
$\ga\in G$ and every 
$$\D (z,z')\in K\cap \Big(Z(s(\ga))\times Z(s(\ga))\Big).$$  Moreover, if we
denote by $K^n$ the set of $(z, z') \in Z*Z$ such that there exists elements
$(z,z_1), (z_1,z_2), \dots (z_{n-1},z')$ in $K$,
we have $\D\cup_{n\geq 0} K^n = Z*Z$ (where ($K^0$ is the diagonal of $Z*Z$).

In order to be able to use the previous corollary \ref{discond-neg2}, we shall
enlarge $Z$ by adding the set of edges. We put $Z^0 = Z$, and  $Z^1 = K/\!\sim$ is
the quotient Borel space of $K$ where we identify an edge with the opposite edge.
The elements of
$Z^1$ are viewed as (non-ordered) subsets of the fibres of $Z$. Let $\tilde Z = Z^0
\sqcup Z^1$ be the disjoint union of $Z^0$ and $Z^1$ with the obvious Borel
structure and
$G$-action. Note that $Z^0$ is the space of {\it vertices} and that $Z^1$ is the
space of {\it unoriented edges}. We endow $\tilde Z$ with a structure of $G$-bundle of metric spaces
in the following way:  we
define the metric $\tilde{d}_x$ on $\tilde{Z}(x)$ by $\tilde{d}_x |_{Z(x)} = d_x$,
 $\tilde{d}_x(z, a) = 1/2$ if $z\in Z(x)$ belongs to $a\in Z^1(x)$ and 
$\tilde{d}_x(z, a) = 0$ otherwise. It is well known that $\tilde{d}_x$ is a
conditionally negative definite kernel on $\tilde{Z}(x)$ (see \cite{wa} for
instance).

An {\it oriented $G$-bundle of trees} is a $G$-bundle of trees $Z$ together
with a Borel subset $K^+$ of $K$, with $K^+ \cap \overline{K}^+ = \emptyset$,
$K = K^+ \cup \overline{K}^+$ and $G K^+ = K^+$, where $\overline{K}^+ :=
\{(z',z) : (z,z')\in K^+\}$ denotes the set of opposite edges of $K^+$.
One also says that $G$ {\it acts on the bundle of trees $Z$ without inversion}.

\begin{thm}[Theorem 4.2, \cite{as}]\label{tree} Let $(G,C)$ be an ergodic
$r$-discrete measured groupoid with property $T$.
\begin{itemize}
\item[(i)] Let $Z$ be a $G$-bundle of trees. Then $Z$ or $Z^1$ has an invariant
section.
\item[(ii)] Every oriented  $G$-bundle of trees has an invariant section. 
\end{itemize}
\end{thm}

\begin{proof} Observe that $\tilde Z$ is a $G$-bundle of complete metric spaces
satisfying the the hypothesis of Corollary \ref{discond-neg2}, and therefore it
has an invariant section. The dichotomy in (i) follows from ergodicity, and from the
$G$-invariance of $Z^0$ and $Z^1$.

When $G$ acts without inversion, every invariant section $x \mapsto a(x)$ of $Z^1$
gives rise to an invariant section of $Z$, by choosing for instance, for all
$x\in G^{(0)}$, the source of the element in $K^+$ representing $a(x)$.
\end{proof}

\begin{defn}\label{treeable} We say that a $r$-discrete Borel groupoid $G$
is {\it treeable} if there exists a structure of Borel $G$-bundle of trees for the
bundle $r :G \ra G^{(0)}$ equipped with the natural left $G$-action.

A measured $r$-discrete groupoid $(G,C)$ is said to be {\it treeable} if there
exists an inessential reduction $(G|_U, C|_U)$ such that the Borel groupoid
$G|_U$ is treeable in the above sense. If $G$ is equipped with such a structure, we
say that $(G,C)$ is a {\it treed measured groupoid}.
\end{defn}

For an ergodic discrete measured equivalence relation this is the usual notion
of treed equivalence relation sudied by Adams and others (see \cite{ad}, \cite{as}).

\begin{cor}[Theorem 1.8, \cite{as}] \label{notreeable} Let $(G,C)$ be an ergodic
$r$-discrete grou\-poid which has property $T$ and is not proper. Then
$(G,C)$ is not treeable. This holds in particular if $(G,C)$ is a properly
ergodic $r$-discrete groupoid with property $T$. 
\end{cor}

\begin{proof} Assume that $(G,C)$ is treeable. By the previous theorem applied with
$r : Z:= G \ra X:= G^{(0)}$, there is a Borel section $\xi$ for $Z$ or for $Z^1$.
In both cases we see that $(G,C)$ is proper. Indeed, for $x\in G^{(0)}$,
if $\xi$ is a section for $r : G \ra G^{(0)}$ we put $\rho^x = \delta_{\xi(x)}$, and
if $\xi : x' \mapsto \{\xi_0(x'), \xi_1(x')\} \subset G^{x'}$ is a section of $Z^1$, we put
$\rho^x = \frac{1}{2}(\delta_{\xi_0(x)} + \delta_{\xi_1(x)})$. Then $\rho$ is an
invariant system of probability measures for $r$, which means that $(G,C)$ is
proper. 

The last observation of the statement comes from the fact that a 
proper ergodic groupoid is not properly ergodic.
\end{proof}

\begin{rem} Let $(G,C)$ be an ergodic
$r$-discrete groupoid which has property $T$.
As an immediate consequence of Corollary \ref{notreeable}, we see that
$(G,C)$ has a structure of oriented $G$-bundle of trees for the natural
left $G$-action on $G$, if and only if (up to an inessential reduction) it is the
coarse equivalence relation on $X:= G^{(0)}$, that is, any two elements of $X$ are
equivalent.

In one direction this observation is obvious: we choose any oriented tree structure
on the countable set $X$, and take the corresponding constant bundle of oriented
trees over $X$. Conversely, let us assume that $(G,C)$ has a structure of oriented
$G$-bundle of trees. According to Theorem \ref{tree} there exists an invariant
Borel section $\xi : x \mapsto \xi(x) \in G^x$. In particular, $(G,C)$ is
proper, and therefore, by Remark \ref{ess-ergo}, we may assume (up to a
conull set) that
$G^{(0)}$ is an orbit.  In addition, since $\ga\xi\circ s(\ga) = \xi\circ r(\ga)$
for every $\ga\in G$, we have $\D\ga = \xi\circ s(\ga)\big(\xi\circ r(\ga)\big)^{-1}$
so that $G$ is an equivalence relation with only one equivalence class.

In case $(G,C)$ only has a structure of $G$-bundle of trees, it is essentially
transitive, but is not necessarily an equivalence relation. For instance the group
$\Z_2$ is treeable in our sense. In fact if $(G,C)$ is treeable, its
isotropy groups are either all trivial, or all equal to $\Z_2$ (up to null sets).
As we have already seen, the case of triviality corresponds to the existence of an
invariant section for $r: G \ra G^{(0)}$. The other possible case is the
existence of an invariant section $x \mapsto \{\xi_0(x), \xi_1(x)\} \subset G^x$.
If there is a non trivial element $\ga$ in the isotropy group $G(x)$ of $x\in
G^{(0)}$, we have $\ga\xi_0(x) = \xi_1(x)$
and therefore $G(x)$ is the group with two elements. By ergodicity, we get that
$G(x)$ is almost everywhere trivial, or almost everywhere equal to $\Z_2$.
\end{rem}

\subsection{Property $T$ and cocycles}

In this last subsection, we study homomorphisms (also called cocycles in case of
group actions) from property $T$ groupoids into locally compact groups
having the Haagerup property.

\begin{defn}\label{Haag-prop} We say that a locally compact group $H$ has the {\it
Haagerup property} (or is {\it a-$T$-menable}) if there exists a continuous function
$\psi$ on
$H$, which is conditionally negative definite and proper, that is, $\lim_{h\ra
\infty}
\psi(h) = +\infty$.
\end{defn}

\begin{exs} For a detailed study of this notion we refer to the book \cite{ccjjv}.
A first list of examples can be found in \cite[Examples 1.2]{ccjjv}. For
instance free groups (and more generally groups acting properly on locally finite
trees, or $\R$-trees), Coxeter groups, amenable groups, have the Haagerup property.
For discrete groups, this property is preserved by amalgamation over a common
finite subgroup \cite[Proposition 6.2.3]{ccjjv}.
\end{exs}

When $G$ is an amenable locally compact group, the following lemma is contained
in Proposition \ref{proper-am}.

\begin{lem}\label{haag-comp} Let $G$ be a locally compact group having the Haagerup
property and let $(X,\mu)$ be an ergodic Kazhdan measured $G$-space. We assume that
for every representation $L$ of $X\croi\, G$ we have $H^1((X\croi\, G, [\mu]), L) = 0$.
Then the $G$-space $(X,\mu)$ is proper. In other words, it is essentially
transitive with compact stabilizer.
\end{lem}

\begin{proof} Due to a result of Varadarajan we can assume that $X$ is a
$\sigma$-compact Hausdorff space and that the action of $G$ is continuous (see
\cite[Theorem 2.19]{zi:ergo} or \cite{ra:topo}). Let $\psi$ be a proper and
continuous  conditionally negative definite 
function on $G$. By \cite[Proposition 14]{hv} there is a representation $\pi$ of $G$
and a $\pi$-cocycle $b$ such that $\psi(g) = \|b(g)\|^2$ for every $g\in G$.
Our assumption implies that the cocycle $(x,g) \mapsto b(g)$, with respect to
the representation $(x,g) \mapsto \pi(g)$, is cohomologous to a coboundary.
Therefore, using Theorem
\ref{triv-cocy}, we get the existence of a Borel subset $E$ of $X$,
of positive measure, such that $\psi$ is bounded on the subset 
$$A := \{g\in G : E \cap gE \not = \emptyset\}.$$ We can take $E$ compact,
so that $A$ is a closed subset of $G$.
Since $\psi$ is proper, observe that $A$ is compact. It follows that the inessential
reduction groupoid $E*G : \{(x,g) : x\in E, g^{-1}x\in E\}$ is contained in $E\times
A$, and has a countably separated space of orbits. Hence, there is a conull orbit in
$E$, and since $[E]$ is conull in $X$, we see that the $G$-space $X$ is essentially
transitive. To end the proof, if suffices to note that the stabilizer of any point
in $E$ is compact.
\end{proof}

\begin{thm}[Theorem 3.2, \cite{jo}] \label{haag} Let $(G,C)$ be an ergodic
$r$-discrete measured groupoid with property $T$, and let $\phi$ be a Borel
homomorphism from $G$ into a locally compact group $H$ having the Haagerup property.
Then $\phi$ is cohomologous to a homomorphism taking values in a compact subgroup of
$H$. 
\end{thm}

\begin{proof} By Theorem \ref{invMackey2} and Lemma \ref{haag-comp},
we get that the Mackey range of the homomorphism $\phi$  is the left action of $H$
on a homogeneous space $H/K$, where $K$ is a compact subgroup of $H$. It follows that
$\phi$ is cohomologous to a homomorphism taking values in $K$.
\end{proof}

\begin{cor} Let $(G,C)$ be an ergodic $r$-discrete Kazhdan groupoid. Then
$r(C)$ contains an invariant measure.
\end{cor}

\begin{proof} Let $\mu\in r(C)$. The Radon-Nikodym derivative $\D \Delta =
\frac{d\mu\circ\lambda} {d(\mu\circ\lambda)^{-1}}$ is a homomorphism from $G$ to
$\R^+$, up to an inessential reduction. By the previous theorem, it is 
cohomologous to the trivial homomorphism, since $\R^+$ has no non-trivial compact
subgroup.
\end{proof}

\begin{rem} Let $G$ be a discrete group containing an infinite subgroup $H$ with
property $T$. Then $G$ admits properly ergodic actions with property $T$. Indeed,
one can easily  construct properly ergodic finite measure preserving actions
of $H$, and it suffices to consider the induced $G$-action.

On the other hand, due to Lemma \ref{haag-comp} an a-$T$-menable group cannot have
property $T$ properly ergodic actions.

Now let us consider a group such as the semi-direct product
$\Z^2 \croi\, SL(2,\Z)$, which is not a-$T$-menable, but does not contain any
infinite property $T$ subgroup. So far, it is unknown whether
$\Z^2 \croi\, SL(2,\Z)$ admits property $T$ properly ergodic actions.
\end{rem}


\begin{thebibliography}{10}

\bibitem{ad} S.~Adams: {\it Trees and amenable equivalence relations}, Erg. Th. Dyn.
Syst., {\bf 10} (1990), 1--14.

\bibitem{as} S.~Adams, R.~Spatzier: {\it Kazdhan groups, cocycles and trees}, Amer.
J. Math., {\bf 112} (1990), 271--287.

\bibitem{dr:amenable}  C.~ Anantharaman-Delaroche and J.~Renault:  Amenable
groupoids, Monographie de l'Enseignement Math\'ematique No {\bf 36},
Gen\`eve, 2000.

\bibitem{ar} W.~Arveson: An invitation to $C^*$-algebras, Graduate Texts in Mathematics n$^o$ {\bf 39},
Springer-Verlag, 1976.

\bibitem{ccjjv} P.-A. Cherix, M. Cowling, P. Jolissaint, P. Julg, A.~Valette:
Groups with the Haagerup property, Gromov's a-$T$-menability, Progress in Mathematics
 no {\bf 197}, Birkh\"{a}user, 2001.

\bibitem{del} P.~Delorme: {\it $1$-cohomologie des repr\'esentations unitaires
des groupes de Lie semi-simples et r\'esolubles. Produits tensoriels continus
et repr\'esentations}, Bull. Soc. Math. France, {\bf 105} (1977), 281--336.

\bibitem{dix} J.~ÊDixmier: Les $C^*$-alg\`bres et leurs repr\'esentations, snd edition, Gauthier-Villars,
1969.

\bibitem{dye1} H.~Dye: {\it On groups of measure transformations I}, Amer. J. Math.,
{\bf 81} (1959), 119--159.

\bibitem{fhm} J.~Feldman, P.~Hahn, C.~C.~Moore: {\it Orbit structure and countable
sections for actions of continuous groups},  Adv. in Math., {\bf 28} (1978),
186--230.

\bibitem{fm1} J.~Feldman, C.~C.~Moore: {\it Ergodic equivalence relations, cohomology
and von Neumann algebras I}, Trans. Amer. Math. Soc., {\bf 234} (1977), 289--324.

\bibitem{fm2} J.~Feldman, C.~C.~Moore: {\it Ergodic equivalence relations, cohomology
and von Neumann algebras II}, Trans. Amer. Math. Soc., {\bf 234} (1977), 325--359.

\bibitem{fe} J.~M.~G.~Fell: {\it Weak containment and induced representations of
groups}, Can. J. Math., {\bf 14} (1962), 237--268.

\bibitem{fu1} A.~Furman: {Gromov's measure equivalence and rigidity of higher rank
lattices}, Ann. of Math., {\bf 150} (1999), 1059--1081.

\bibitem{fu2} A.~Furman: {Orbit equivalence rigidity}, Ann. of Math., {\bf 150}
(1999), 1083--1108.

\bibitem{ga} D.~Gaboriau: {\it Invariants $\ell^2$ de relations d'\'equivalences et
de groupes}, Publ. Math. I.H.E.S., to appear.

\bibitem{gui:coho} A.~Guichardet: {\it Etude de la $1$-cohomologie
et de la topologie du dual pour les  groupes de Lie \`a radical ab\'elien},
Math. Ann., {\bf 228} (1977), 215--232.

\bibitem{ha1}  P.~Hahn: {\it Haar measure for measure groupoids}, Trans. Amer.
Math. Soc., {\bf 242} (1978), 1--33.

\bibitem{hv} P.~de la Harpe, A.~Valette:  La propri\'et\'e (T) de Kazdhan por les groupes
localement compacts (avec un appendice de M. Burger), Ast\'erisque n$^0$ {\bf 175}, 1989.

\bibitem{jo} P.~Jolissaint: {\it Borel cocycles, approximation properties and
relative property $T$}, Ergod. Th. and Dynam. Syst., {\bf 20} 2000, 483--499.

\bibitem{ka} D.~Kazhdan: {\it Connection of the dual space of a group with the structure of its
closed subgroups}, Funct. Anal. and its Appl., {\bf 1} (1967), 63--65.

\bibitem{ko} G.~K\"{o}the: Topological Vector Spaces I, Springer-Verlag, 1969.

\bibitem{mo:extIII} C.~C.~Moore: {\it Group extensions and cohomology for locally
compact groups. III},  Trans. Amer. Math. Soc., {\bf 221} (1976), 34--58.

\bibitem{mo:ergo} C.~C.~Moore: {\it Ergodic theory and von Neumann algebras},
in Operator and applications, Part 2 (Kingston, Ont. 1980), 179--226. Amer. Math.
 Soc., 1982.

\bibitem{mu} P.~Muhly: Coordinates in Operator Algebras, in preparation. 

\bibitem{ne} A.~Nevo: {\it Amenable actions and actions with property $T$},
M. Sc. Thesis, Hebrew University, 1987 (in Hebrew).

\bibitem{po} S.~Popa: {\it On a class of type $II_1$ factors with Betti numbers
invariants}, Preprint, 2002.

\bibitem{ra:virt} A. Ramsay: {\it Virtual groups and group actions}, Adv. in Math., {\bf 6}
 (1971), 253--322.

\bibitem{ra:nontran} A. Ramsay: {\it Nontransitive quasiorbits in Mackey's analysis
of group extensions}, Acta Math., {\bf 137} (1976), 1--48.

\bibitem{ra:topo} A. Ramsay: {\it Topologies for measured groupoids}, J. Funct. Anal., {\bf 47}
(1982), 314-343.

\bibitem{sch1} K.~Schmidt: Cocycles on ergodic transformation groups,
MacMillan (Company of India Ldt., Dehli), 1977.

\bibitem{sch3} K.~Schmidt: {\it Amenability, Kazhdan property $T$, strong
ergodicity and invariant means for ergodic group-actions}, Ergod. Th. and
Dynam. Sys., {\bf 1} (1987), 223--236.

\bibitem{sch2} K.~Schmidt: Algebraic ideas in ergodic theory, C. B. M. S.,
Regional Conf. Series in Math. 
no {bf 76}, Amer. Math. Soc., 1990.

\bibitem{wa} Y.~Watatani: {\it Property $(T)$ of Kazhdan implies property $(FA)$
of Serre}, Math. Japonica, {\bf 27} (1981), 97--103.

\bibitem{we} J.~Westman: {\it Cohomology of ergodic groupoids}, Trans. Amer. Math.
Soc., {\bf 146} (1969), 465--471.

\bibitem{zi:amen} R.J. Zimmer: {\it Amenable ergodic group actions and an
application to Poisson boundaries of random walks}, J. Funct. Anal., {\bf 27}
(1978), 350--372.

\bibitem{zi:induce} R.J. Zimmer: {Induced and amenable ergodic actions of Lie groups},
Ann. Scien. Ec. Norm. Sup., {\bf 11} (1978), 407--428.

\bibitem{zi:ergo} R.J. Zimmer: Ergodic theory and semisimple groups, Birkh\"{a}user,
1984.

\bibitem{zi:cohom} R.J.~Zimmer: {\it On the cohomology of  ergodic actions of
semisimple Lie groups and discrete subgroups}, Amer. J. Math., {\bf 103} (1981),
937--951.

\bibitem{zi:kazhdan} R.J.~Zimmer: {\it Kazhdan groups acting on compact manifolds},
Invent. Math., {\bf 75} (1984), 425--436.

\end{thebibliography}
\end{document}